\documentclass[12pt]{amsart}
\usepackage[utf8]{inputenc}
\usepackage{graphicx}
\usepackage{subcaption}
\usepackage{amsmath}
\usepackage{amssymb}
\usepackage{upgreek}
\usepackage{mathtools}
\usepackage{centernot}
\usepackage{stmaryrd}
\usepackage{microtype}
\usepackage{esvect}
\usepackage{amsthm}
\usepackage{upgreek}
\usepackage{comment}
\usepackage{amsfonts}
\usepackage{physics}
\usepackage{tikz-cd}
\usepackage[utf8]{inputenc}
\usepackage[english]{babel}
\usepackage{graphicx}
\usepackage{mathtools}
\usepackage{latexsym}
\usepackage{dsfont}
\usepackage{enumerate}
\usepackage[shortlabels]{enumitem}
\usepackage[hidelinks]{hyperref}
\usepackage{multirow}
\usepackage{url}
\usepackage{setspace}

\usepackage{xargs}                  
\usepackage{xcolor}
\usepackage{times}

\usepackage[colorinlistoftodos,prependcaption,textsize=tiny]{todonotes}

\newtheorem{theorem}{Theorem}[section]

\theoremstyle{definition} 

\newtheorem{conjecture}[theorem]{Conjecture}
\numberwithin{equation}{subsection}

\makeatletter
\newsavebox{\@brx}
\newcommand{\llangle}[1][]{\savebox{\@brx}{\(\m@th{#1\langle}\)}%
  \mathopen{\copy\@brx\kern-0.5\wd\@brx\usebox{\@brx}}}
\newcommand{\rrangle}[1][]{\savebox{\@brx}{\(\m@th{#1\rangle}\)}%
  \mathclose{\copy\@brx\kern-0.5\wd\@brx\usebox{\@brx}}}
\makeatother

\usepackage{subfiles}
\usepackage[labelformat=simple]{subcaption}

\makeatletter
\def\namedlabel#1#2{\begingroup
 #2%
 \def\@currentlabel{#2}%
 \phantomsection\label{#1}\endgroup
}
\makeatother

\title[Totally Geodesic Spanning Surfaces]{Totally Geodesic Spanning Surfaces of Knots and Links in 3-Manifolds}

\author[Benjamin Shapiro]{Benjamin Shapiro}

\begin{document}

\thispagestyle{empty}

\maketitle

\begin{abstract}
    I construct infinite families of knots and links with totally geodesic spanning surfaces, which we call \emph{TGS knots} and \emph{TGS links}, in various 3-manifolds. These 3-manifolds include thickened orientable surfaces, the sphere cross the circle, lens spaces, and the solid torus. The totally geodesic spanning surfaces of knots embedded in thickened orientable surfaces are also examples of totally geodesic spanning surfaces of virtual knots.
\end{abstract}

\section{Introduction}\label{intro}

\subsection{Motivation and Purpose}

A knot $K$ in a 3-manifold manifold $M$ is an embedding of the circle $S^1$ in $M$. A link $L$ is an embedding of some finite number of copies of $S^1$ in some 3-manifold $M$. In this paper, we will look at embeddings of knots and links in a collection of 3-manifolds, including thickened surfaces, $S^2 \times S^1$, lens spaces, and the solid torus.

A spanning surface $S$ of a knot or link $L$ in a 3-manifold $M$ is a surface in $M$ such that the boundary of $S$ is equal to $L$. That is, $S \cap L = L = \partial S$. A link $L$ is said to be hyperbolic in a 3-manifold $M$ if its complement $M \setminus L$ is a hyperbolic 3-manifold. A surface $S$ embedded in a hyperbolic 3-manifold $M$ is said to be totally geodesic if it is isotopic to a surface that lifts to a set of geodesic planes in hyperbolic 3-space, denoted $H^3$.

Totally geodesic surfaces are interesting surfaces to study due to their cleanliness and beauty. These surfaces are able to elicit a 2-dimensional hyperbolic geometry from the 3-dimensional hyperbolic geometry of $H^3$. Their cleanliness comes from the fact that geodesic planes in $H^3$ correspond to vertical planes and perfect hemispheres in the upper-half space model of $H^3$, and such simplicity in $H^3$ is not common. Totally geodesic surfaces are important since they play the same role in 3-manifolds as geodesics play in 2-manifolds. Geodesics denote line segments of minimal length between two points along a surface, and totally geodesic surfaces denote surfaces of minimal area. Certain totally geodesic surfaces, specifically thrice-punctured spheres, are also very important when it comes to performing surgeries on a hyperbolic 3-manifold. We are able to cut a hyperbolic 3-manifold along a thrice-punctured sphere, glue the thrice-punctured sphere back to itself in some non-identity mapping, and the structure of the totally geodesic surface will preserve geometry in other sections of the manifold.

Given this interest in totally geodesic surfaces, we want to figure out more about when and where we can find them in various 3-manifolds. Note, Myers' work in \cite{Myers} and Adams' work in \cite{Augment} and \cite{thrice} can be combined to show that totally geodesic surfaces with no other constraints can be constructed quite easily. In \cite{Myers}, Myers proves every closed orientable 3-manifold contains a hyperbolic knot. Then, Adams proves in \cite{Augment} that we can always find a second component which bounds a twice-punctured disk. Note, a twice-punctured disk is equivalent to a thrice-punctured sphere. Lastly, Adams proves in \cite{thrice} that thrice-punctured spheres are always totally geodesic in hyperbolic 3-manifolds. Hence, we know that every 3-manifold contains at least one link with at least one totally geodesic surface, and given any hyperbolic knot we can add a single component in order to yield a totally geodesic surface.

Therefore, in this paper we are specifically interested in totally geodesic spanning surfaces of knots and links. Define a \emph{TGS knot} (or \emph{TGS link}) to be a knot $K$ embedded in some 3-manifold $M$ such that $M \setminus K$ is hyperbolic and $K$ has at least one totally geodesic spanning surface. We want to know how common TGS knots and links are and in which 3-manifolds we can find them, and we conjecture they are present in many places.

\begin{conjecture}
Every closed orientable 3-manifold contains at least one TGS link.
\end{conjecture}

To begin building towards this large conjecture, in this paper we will look at specific 3-manifolds and prove the existence of infinite families of TGS knots and links in those manifolds. The purpose of ensuring these families of knots and links are infinite is to show TGS knots and links are more common than simply finding  only a few special cases where totally geodesic spanning surfaces exist.

In \cite{KaplanKelly}, Kaplan-Kelly gives examples of TGS links in thickened surfaces which contain two totally geodesic spanning surfaces. She does so by showing the complements of these links can be decomposed into two right-angled generalized polyhedra glued together by rotation on their faces. In this paper, we will construct a different set of TGS knots in thickened surfaces. We will use rigid orbifolds to prove the spanning surfaces of these knots are totally geodesic, rather than right-angled generalized polyhedra compositions. 

Furthermore, the authors of \cite{CarterKamadaSaito}, \cite{KamadaKamada}, and \cite{Kauffman} prove that any virtual knot can be realized as a classical knot in a thickened surface. Virtual knots are an extension of classical knots which allow for a third option at each crossing, called a virtual crossing. Therefore, the TGS knots in thickened surfaces and their totally geodesic spanning surfaces discussed in this paper are examples of totally geodesic spanning surfaces of virtual knots. These examples follow after Kaplan-Kelly's example of a knot in a thickened surface with two totally geodesic spanning surfaces in \cite{KaplanKelly}.

\subsection{Theorems and Methods}

Now, we will introduce some important theorems and how we will use them to show certain knots and links are hyperbolic in a given 3-manifold and to show certain surfaces are totally geodesic.

\begin{theorem}\label{small2018}
(Adams, Albors-Riera, Haddock, Li, Nishida, Reinoso, and Wang \cite{thickened}) Let $M$ be a finite volume hyperbolic 3-manifold, possibly with cusps such that any
boundary is total geodesic. Let $S$ be an essential, closed surface in $M$ with neighborhood $N$. For
any link $L$ that is prime in $N$ with a cellular alternating projection on $S$, the manifold $M \setminus L$ is
hyperbolic.
\end{theorem}

We will use Theorem \ref{small2018} to help construct knots in thickened surfaces which have hyperbolic complements. To do so, we will need to ensure the knot is prime and cellular alternating. A projection of a knot $K$ from a thickened surface $S \times I$ to $S$ is said to be fully alternating if it is alternating on $S$ and the interior of every complementary region is an open disk.

\begin{theorem}\label{HowiePurcell}
 (Howie and Purcell \cite{hypLinks}) Given a projection surface $F$ in a compact, orientable, irreducible 3-manifold $Y$ and a link $L$ in $Y$. If $F$ has genus at least one, and the regions in the complement of $\pi(L)$ on $F$ are disks, and the representativity $r(\pi(L), F) > 4$, then $Y \setminus L$ is hyperbolic.
\end{theorem}

Representativity is the minimum number of intersections any curve which bounds a compressing disk has with the link $L$ in the projection surface. We will use Theorem \ref{HowiePurcell} and the notion of representativity to construct knots in other 3-manifolds which are hyperbolic. In this paper, we will always have the projection surface $F$ be a torus.

\begin{theorem}\label{Thurston}
(Thurston \cite{dehn}) Let $M(u_{1},u_{2},\dots ,u_{n})$ denote the manifold obtained from $M$ by filling in the $i$-th boundary torus with a solid torus using the slope $u_{i}=p_{i}/q_{i}$ where each pair $p_{i}$ and $q_{i}$ are coprime integers. Then,
$M(u_{1},u_{2},\dots ,u_{n})$ is hyperbolic as long as a finite set of exceptional slopes is avoided for the $i$-th cusp for each $i$.
\end{theorem}

Theorem \ref{Thurston} is extremely powerful. In this paper, we only use it to prove any link in which is hyperbolic in the torus is not hyperbolic in finitely many lens spaces. Therefore, it will be hyperbolic in infinitely many lens spaces.

Formally, an orbifold is defined to be a space locally modelled on $\mathbb{R}^n$ modulo finite group actions. In this paper, we will use orbifolds to study 3-manifolds quotiented by a subset of their symmetries. An orbifold is said to be rigid if there exists only one possible hyperbolic metric that generates the orbifold.

\begin{theorem}\label{AdamsSchoenfeld}
(Adams and Schoenfeld \cite{orbifolds}) Let M be a closed orientable 3-manifold, covering a 3-orbifold O such that M contains a link L covering a collection J of simple closed curves and arcs beginning and ending on the order two singular set of O, and such that $O \setminus J$ is a hyperbolic 3-orbifold O'. Then if O' contains a rigid hyperbolic 2-orbifold Q for which each simple closed component of J is a single puncture (singularity with $m_i$ = $\infty$) and each arc corresponds to a corner reflector with $n_j$ = $\infty$, then Q is isotopic in O' to a 2-orbifold that lifts to a totally geodesic Seifert surface for the link L in M.
\end{theorem}

We will use Theorem \ref{AdamsSchoenfeld} to show that various spanning surfaces are totally geodesic. Given a spanning surface, we will take the quotient of the surface by some of its symmetries. This will yield some orbifold. If we can then show the orbifold is rigid, then Theorem \ref{AdamsSchoenfeld} tells us the orbifold lifts to a totally geodesic surface, with that surface being the spanning surface we began with. The technique for proving rigidness used in \cite{orbifolds} which we will also use in this paper involves using the orders of the orbifold's corner reflectors and calculating the sum of their reciprocals. A corner reflector is an intersection point between two symmetry axes or a symmetry axis and the boundary of the surface. Its order is dictated by the orders of the symmetries which produced that corner reflector.

\begin{theorem}\label{AdamsThrice}
(Adams \cite{thrice}) Let M be a compact orientable 3-manifold such that M is a finite volume hyperbolic 3-manifold. Let S be an incompressible thrice-punctured sphere properly embedded in M. Then S is isotopic to a thrice-punctured sphere S' properly embedded in M such that int(S') is totally geodesic in the hyperbolic structures on int(M).
\end{theorem}

Theorem \ref{AdamsThrice} tells us any surface which is topologically a thrice-punctured disk must be totally geodesic.

\begin{theorem}\label{KojimaGeodisic}
(Kojima \cite{Kojima}) The underlying space of a nonsingular part N of an orientable hyperbolic 3-cone-manifold C of finite volume carries a complete negatively curved metric. In particular it is homeomorphic to an interior of a compact irreducible atoroidal 3-manifold with toral boundary which admits no Seifert fibrations. Moreover, it admits a complete hyperbolic structure of finite volume. 
\end{theorem}

Given a hyperbolic link, we will use Theorem \ref{KojimaGeodisic} to prove that adding another component which is a geodesic in its complement to that link yields another hyperbolic link.

In the following sections, I construct families of knots and links embedded in various 3-manifolds and then apply these theorems to prove they are hyperbolic and that specific spanning surfaces of these knots and links are totally geodesic.

\subsection{Acknowledgements}

I would like to thank my advisor, Colin Adams, for his valuable guidance and insight.

\section{TGS Knots in Thickened Orientable Surfaces}\label{thick}

In this section, we will look at infinite families of TGS knots which are embedded in thickened, orientable surfaces of arbitrary genus. 

\subsection{A First Infinite Family in Each Thickened Surface}\label{th1}

\begin{figure}[htbp]
\begin{center}
\includegraphics[width=4cm]{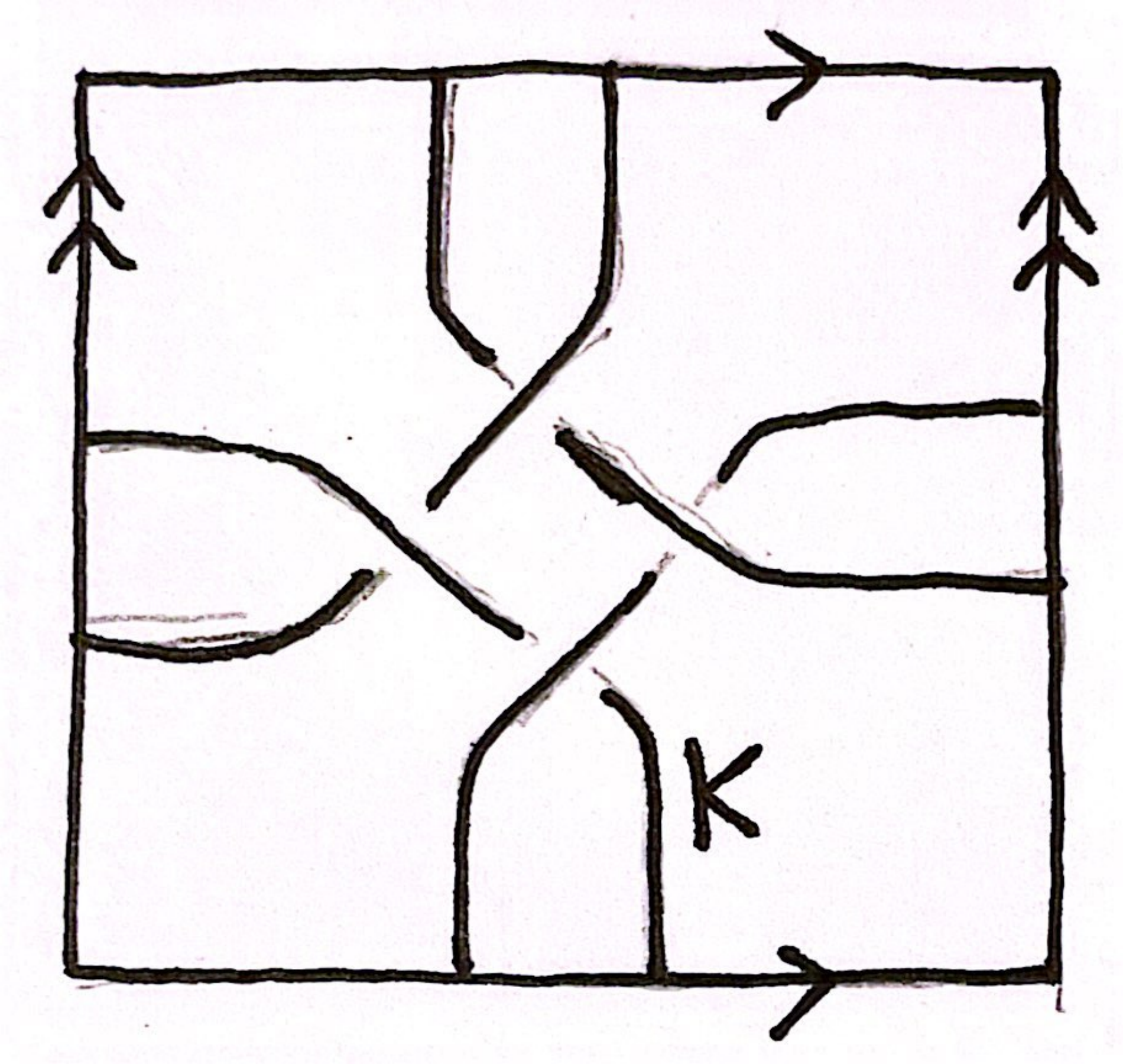}

\caption{}
\label{ThTorSq}
\end{center}
\end{figure}

We will begin by constructing a TGS knot embedded in a thickened torus. Let $G$ be a square gluing diagram with opposite edges identified which represents a torus $F$. We will consider the thickened surface $F \times [-1, 1]$ with $G = F \times \{0\}$. Next, we will construct a knot projection $K$ on $G$. Begin by forming strands into a square in the center of $G$ such that the crossings are alternating. Then, extend the two strands at each crossing to the corresponding edge of $G$, as shown in Figure \ref{ThTorSq}. We will refer to these strands that extend to the edges of $G$ as ``arms" of $K$.

\begin{figure}[htbp]
\begin{center}
\includegraphics[width=4.5cm]{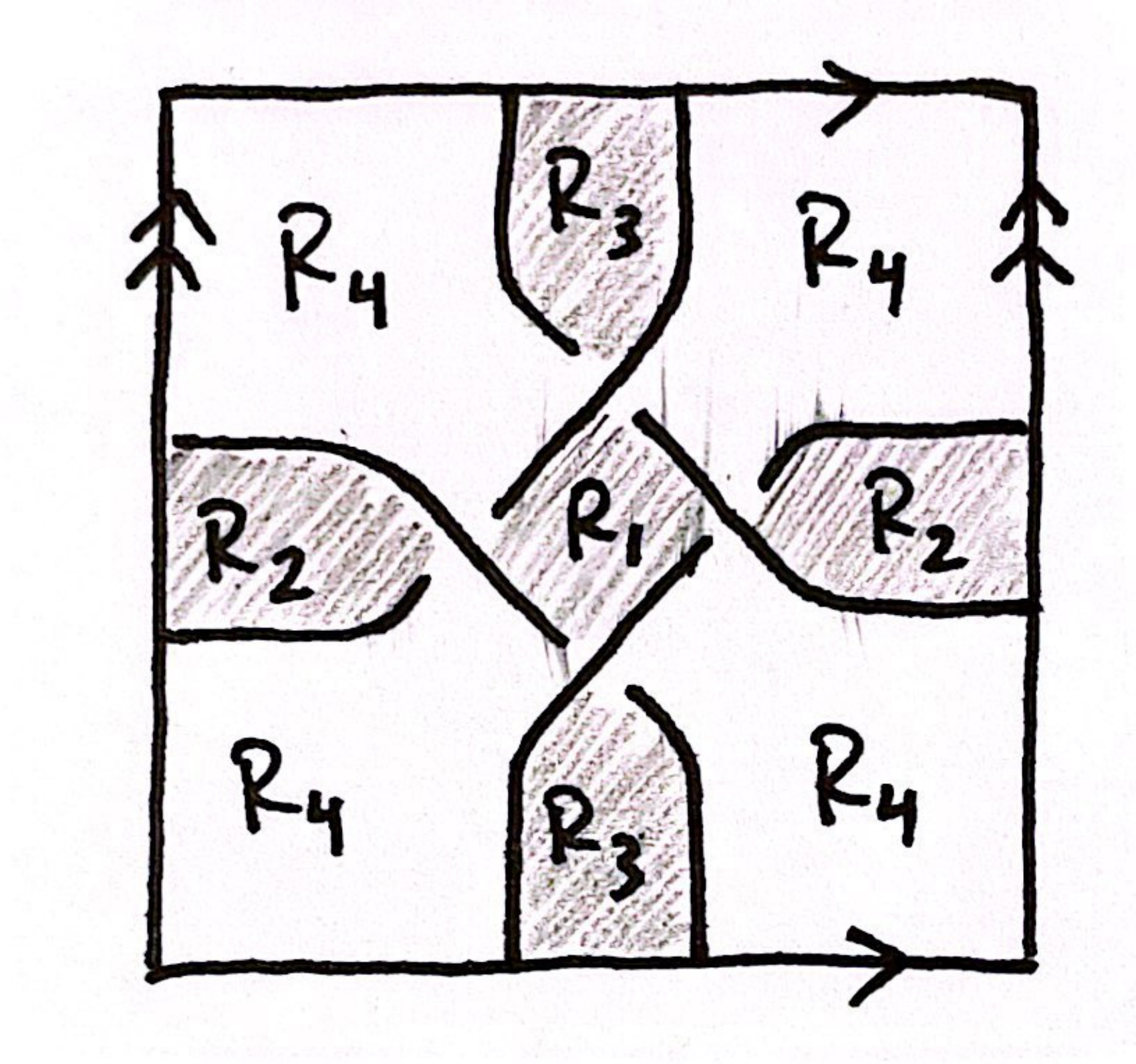}

\caption{}
\label{ThTorSqRegs}
\end{center}
\end{figure}

We claim that $K$ is a hyperbolic knot and $K$ contains a totally geodesic surface in its complement. Note, $K$ is alternating. The complement of $K$ is composed of four regions $R_1, R_2, R_3, R_4$ as in Figure \ref{ThTorSqRegs}. $R_1$ is a square, $R_2$ and $R_3$ are bigons, and $R_4$ is an octagon. Therefore, the interior each of these four regions is a topological open disk. Hence, $K$ is cellular alternating.

\begin{figure}[htbp]
\begin{center}
\includegraphics[width=9cm]{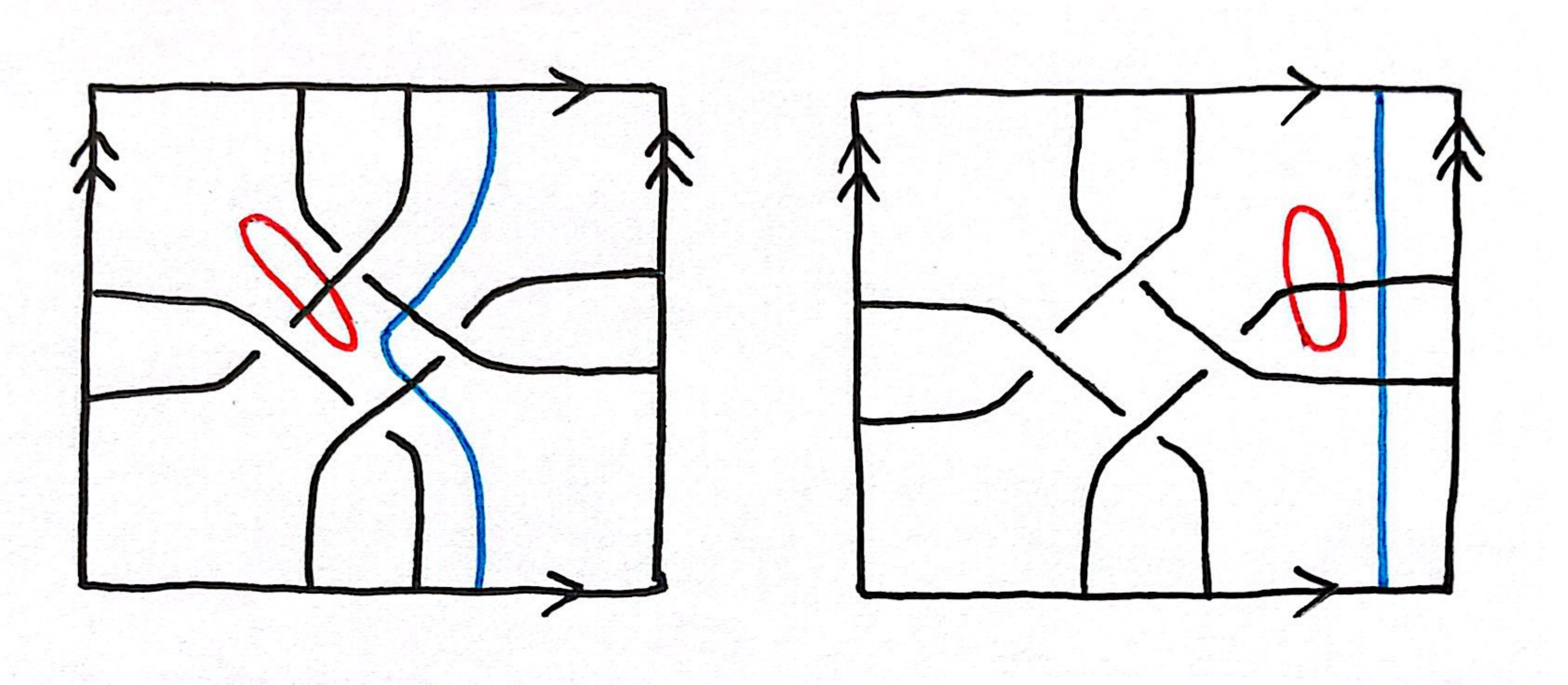}

\caption{}
\label{ThTorSqObvPr}
\end{center}
\end{figure}

Next, we will show that $K$ is obviously prime. In \cite{thickened} the authors define obviously prime as follows: ``A reduced, fully alternating projection $P$ of a link $L$ in $S \times I$ onto $S$ is obviously prime if every disk in the projection surface with boundary intersecting $P$ transversely at two points intersects P in an embedded arc." In our torus gluing diagram $G$, a simple closed curve which crosses over an edge of $G$ once does not bound a disk since it includes either a meridian or a longitude of the torus. Let $U$ be a simple closed curve which intersects $K$ exactly two times. Note, since the four complementary regions of $K$ are topological open disks, $U$ can only pass through two adjacent regions. By adjacent here we mean they meet along a strand of the knot, not along a crossing. Additionally, $R_1, R_2,$ and $R_3$ are all only adjacent to $R_4$. So, $U$ must pass through $R_4$ and one other region. Now, we can split this into two cases. In the first case, $U$ runs through the central square region $R_1$. In the second, $U$ runs through a bigon region, either $R_2$ or $R_3$. Let's consider the first case. We will begin drawing $U$ in the interior of $R_1$. We must cross some strand on the boundary of $R_1$ to pass into $R_4$. Then, to close $U$ we can either cross over that same strand or pass through an edge of $G$ to cross over some other strand on the boundary of $R_1$. These two options are depicted in red and blue in the left diagram of Figure \ref{ThTorSqObvPr}, respectively. The red curve intersects $K$ in an embedded arc and the blue curve does not bound a disk. Therefore, the conditions to be obviously prime are met in this case. Next, let's consider the second case, which uses similar reasoning. Without loss of generality, let's assume $U$ passes through $R_2$. In the this case, $U$ can either intersect one of the boundary strands of $R_2$ twice, or both boundary strands of $R_2$ once each. These two options are depicted in red and blue in the right diagram of Figure \ref{ThTorSqObvPr}, respectively. Again, the red curve intersects $K$ in an embedded arc and the blue curve does not bound a disk. Hence, the conditions are met in this case as well. It follows that $K$ is obviously prime. 

In \cite{thickened}, the authors also prove that a cellular alternating knot in a thickened surface is prime if and only if any reduced cellular alternating projection of it is obviously prime. Hence, the knot of which $K$ is a projection is prime. It follows by Theorem \ref{small2018} that this knot is hyperbolic in $F \times [-1,1]$.

\begin{figure}[htbp]
\begin{center}
\includegraphics[width=9cm]{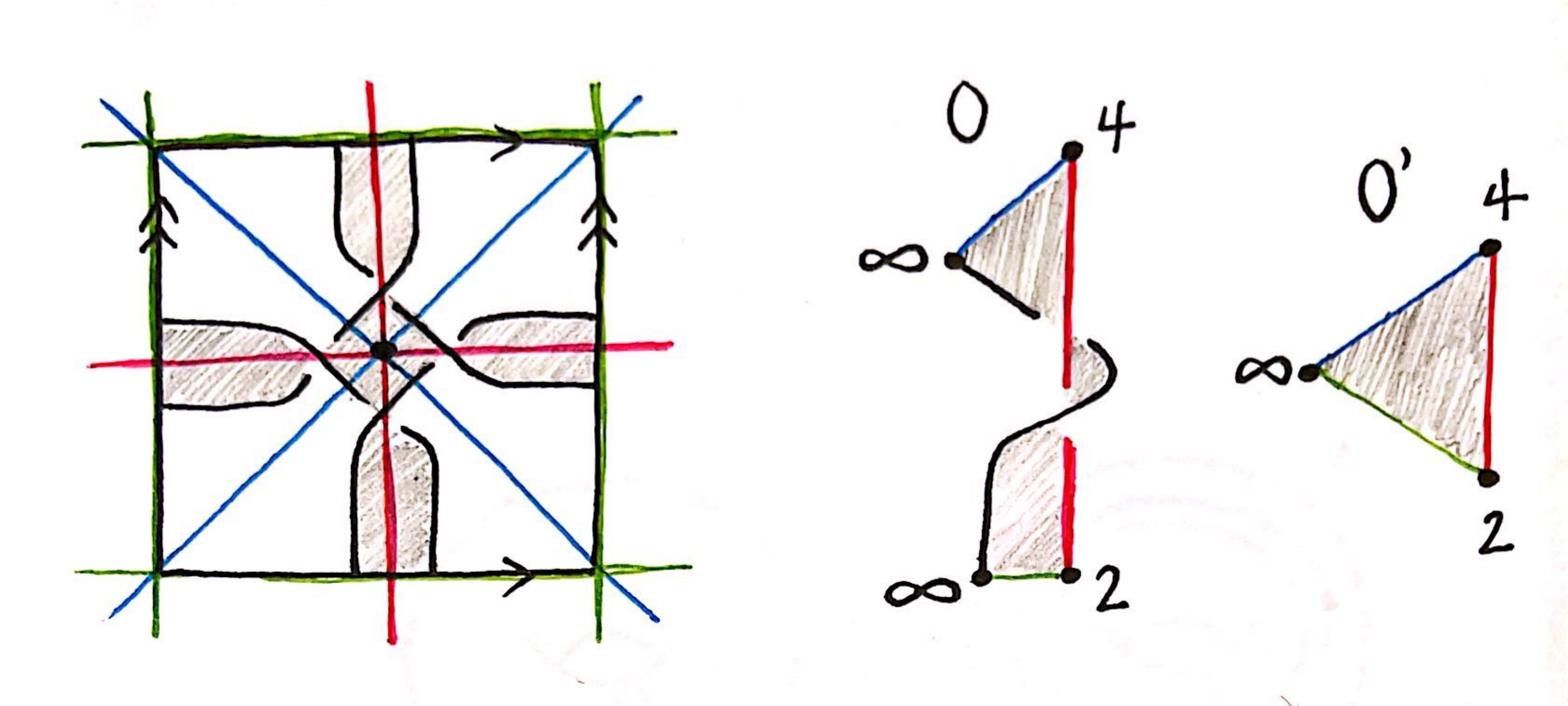}

\caption{}
\label{ThTorSqOrb}
\end{center}
\end{figure}

Now, we must show that $K$ has a totally geodesic spanning surface. Consider the surface $S = R_1 \cup R_2 \cup R_3$. $S$ possesses rotational symmetry of degree 4 about the center point of $R_1$. $S$ also possesses reflective symmetry over the three lines, call them $L_1, L_2, L_3$ shown in the left diagram of Figure \ref{ThTorSqOrb}; $L_1$ runs vertically through one bigon and through opposite crossings of the central square, $L_2$ runs diagonally through opposite edges of the central square, and $L_3$ runs horizontally through one bigon. Quotienting $S$ by these three symmetries gives an orbifold $O$ as in the center diagram of Figure \ref{ThTorSqOrb}. $O$ has four corner reflectors with orders 4, 2, $\infty$, and $\infty$. In this figure, the black knot segment and $S$ both wrap around $L_1$. This is due to the crossing through which $L_1$ runs in the original diagram of $K$. Thus, $O$ can be ``untwisted" to give an isotopic orbifold that lies flat in the plane. Since $O$ is embedded in the complement of $K$, the knot segment of $O$ is an open boundary of $O$. Therefore, the two corner reflectors of degree $\infty$ can be reduced along the segment of $K$ to a single point. This reduction gives an orbifold $O'$ equivalent to $O$, as shown in the right diagram of Figure \ref{ThTorSqOrb}, with three corner reflectors. Note, $\frac{1}{4} + \frac{1}{2} + \frac{1}{\infty} = \frac{3}{4} < 1$. Therefore, by Theorem \ref{AdamsSchoenfeld}, $O'$ is a rigid hyperbolic 2-orbifold which lifts to $S$ in $(F \times [-1,1]) \setminus K$. It follows that $S$ is totally geodesic, and thus $K$ is a TGS knot.

\begin{figure}[htbp]
\begin{center}
\includegraphics[width=4.5cm]{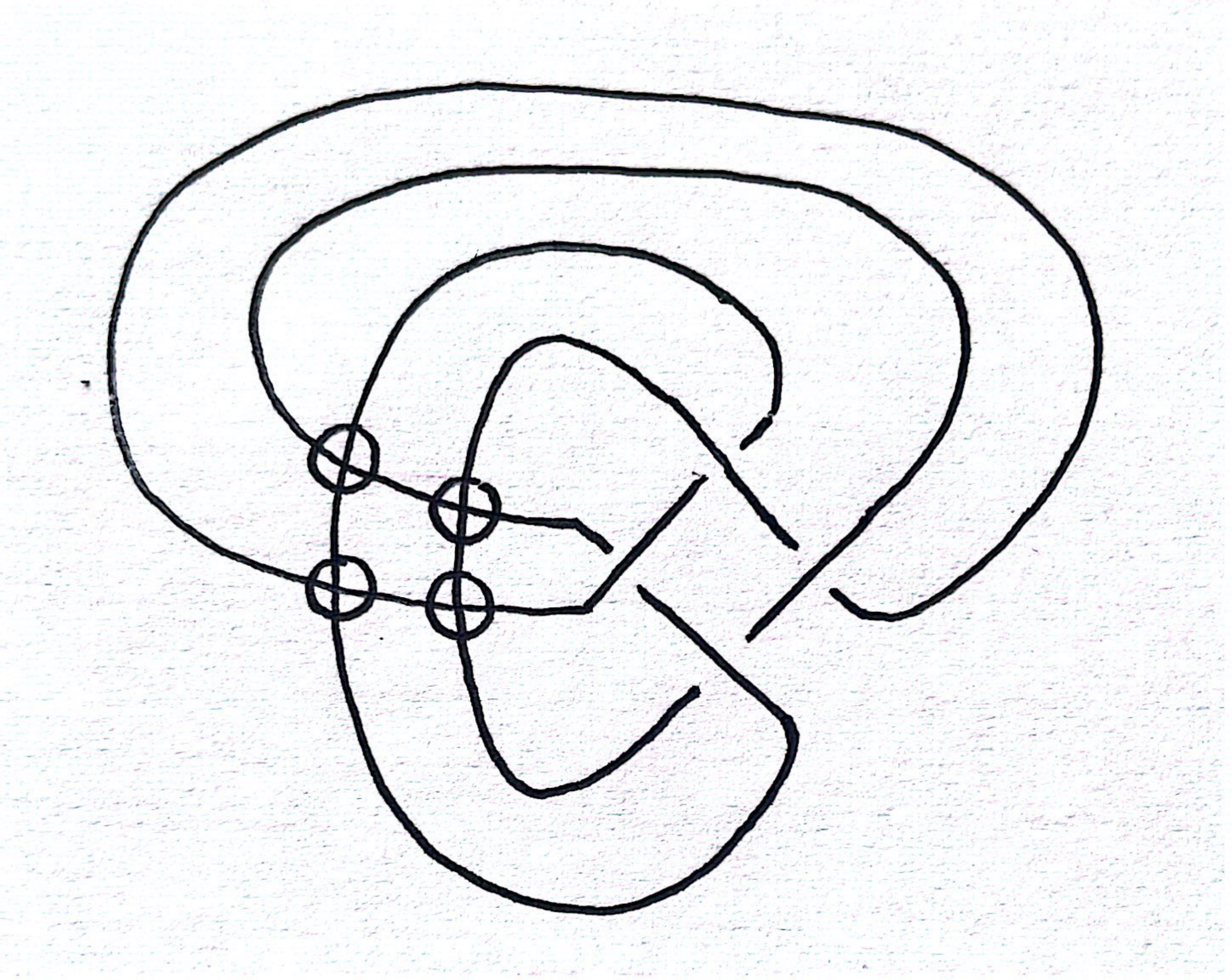} \includegraphics[width=4.5cm]{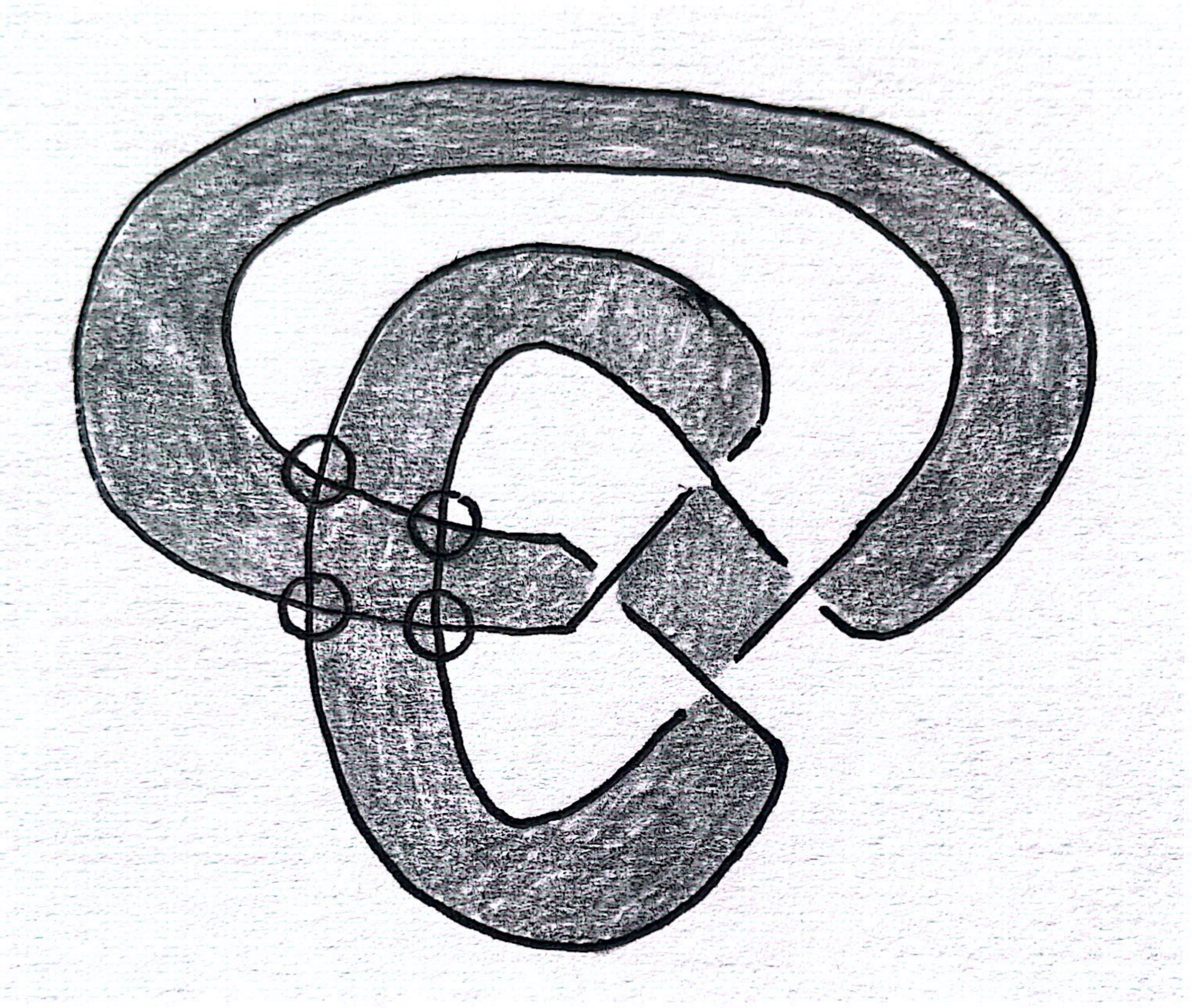}

\caption{}
\label{virtual}
\end{center}
\end{figure}

As we mentioned in the Section \ref{intro}, the knot $K$ and its spanning surface $S$ correspond to a virtual knot $K_V$ and one of its spanning surfaces, respectively. The left diagram of Figure \ref{virtual} depicts $K_V$. Since $S$ is totally geodesic, the spanning surface $S_V$ of $K_V$ shown in the right diagram of Figure \ref{virtual} is also a totally geodesic spanning surface. Thus, we have the first example of a totally geodesic spanning surface of a virtual knot.

\begin{figure}[htbp]
\begin{center}
\includegraphics[width=4cm]{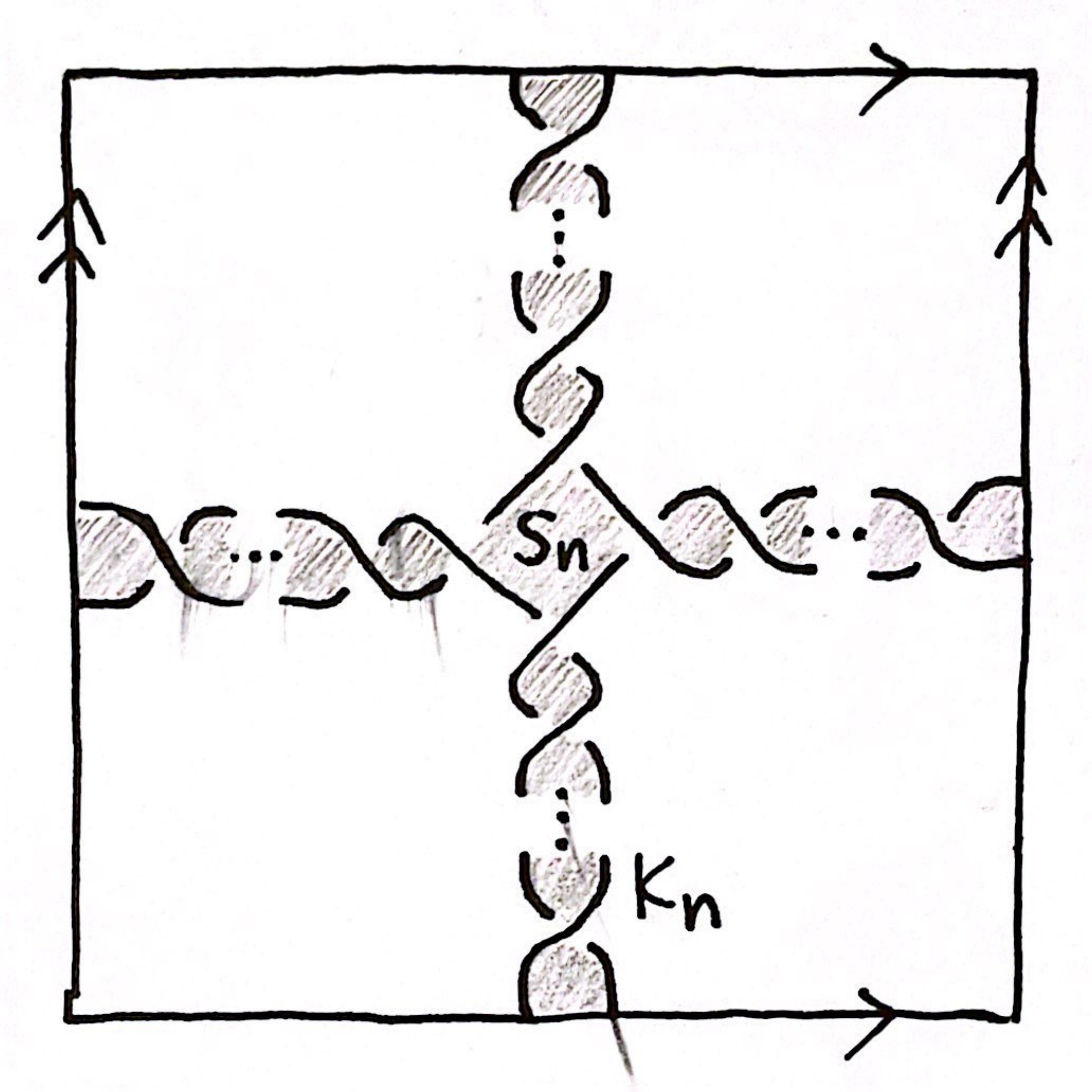}

\caption{}
\label{ThTorSqKn}
\end{center}
\end{figure}

Next, we are able to construct an infinite family of TGS knots in $F \times [-1,1]$ based on the example above. Take $K$ and add $n$ half-twists to each arm, ensuring the resulting knot remains alternating. Call this new knot $K_n$. Then, let $S_n$ be the spanning surface of $K_n$  which includes the interiors of the central square and each arm as in Figure \ref{ThTorSqKn}.

With the added twists, each arm consists of $n+1$ bigons. Each complementary region is still a disk, so $K_n$ is cellular alternating. Additionally, all of the bigons which form the arms in addition to the central square are only adjacent to the outer region. Therefore, the same argument holds as before which proves $K_n$ is obviously prime. That is, any simple closed curve which intersects $K$ twice runs through the outer region and either the central square region or one bigon region. In both cases, the curve either bounds a disk which intersects $K_n$ along an embedded arc or does not bound a disk. Therefore, $K_n$ is obviously prime. It follows by Theorem \ref{small2018} that $\mathcal{K}_n$, the knot of which $K_n$ is a projection, is hyperbolic in $F \times [-1,1]$.
 
\begin{figure}[htbp]
\begin{center}
\includegraphics[width=9cm]{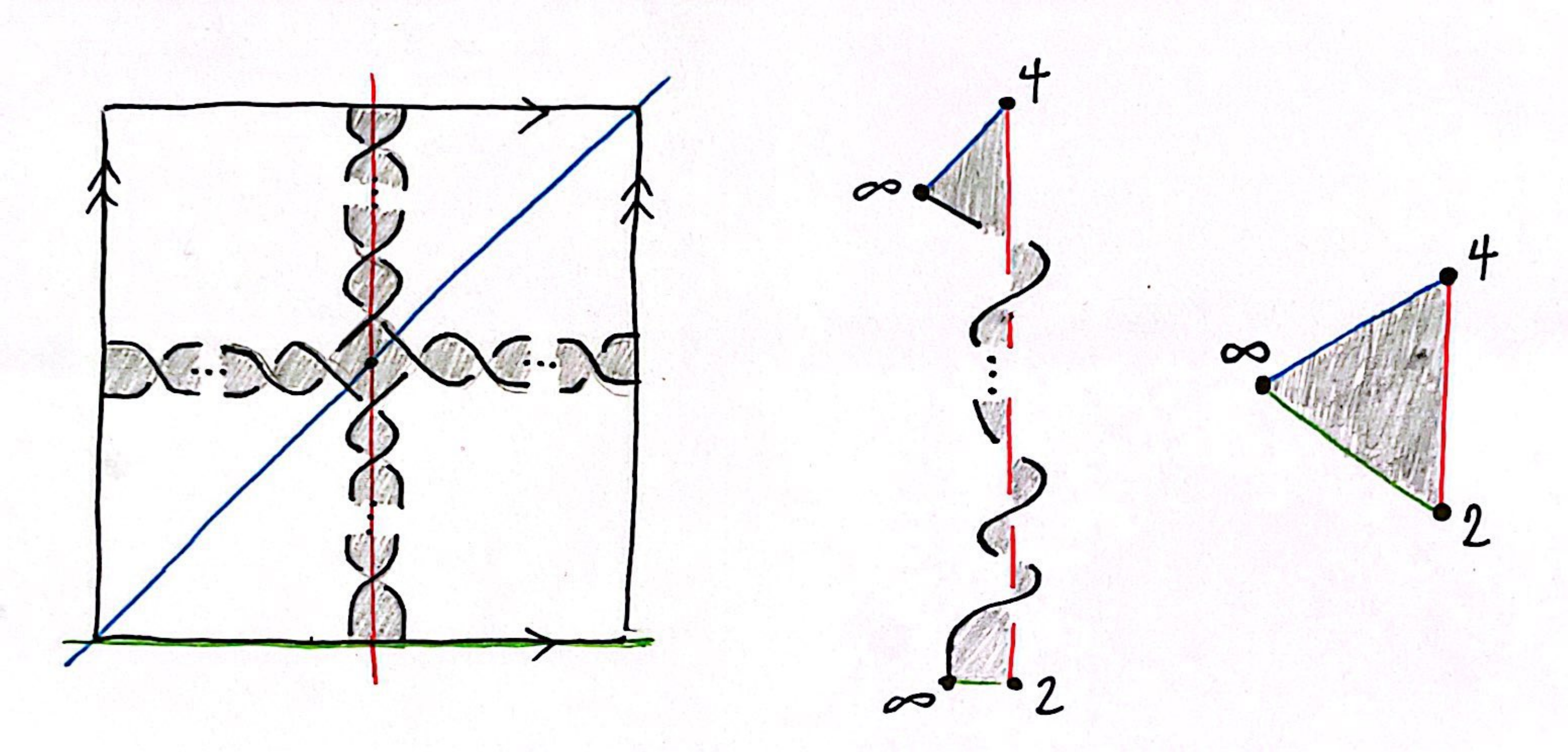}

\caption{}
\label{ThTorSqKnOrb}
\end{center}
\end{figure}

 $S_n$ possesses the same three lines of reflective symmetry $L_1, L_2, L_3$ as $S$ did above, which are shown in the left diagram of Figure \ref{ThTorSqKnOrb}. Then, taking the quotient of $S_n$ by these three lines of reflective symmetry gives the orbifold $O_n$ shown in the center diagram of Figure \ref{ThTorSqKnOrb}. Again, the knot component and $S_n$ wrap around $L_1$ once per crossing of $K_n$ that $L_1$ runs through. Now, we can also perform isotopy moves and untwist $O_n$ so that it lies flat in the plane to show $O_n$ is isotopic to $O$. Therefore, $O_n$ is a rigid orbifold since $O$ is a rigid orbifold. Therefore, $S_n$ must be totally geodesic. It also follows that $K_n$ is a TGS knot. This gives a countably infinite family of TGS knots embedded in the thickened torus.

Next, we will show that this construction can be used to find an infinite family of TGS knots in any thickened genus $g$ surface with $g \geq 2$. Let $G^g$ be a gluing diagram that is a regular $4g$-gon with opposite edges identified by translation. This represents an orientable genus $g$ surface $F^g$. Like before, we will consider the thickened surface $F^g \times [-1, 1]$ with $G^g = F^g \times \{0\}$. Similarly to the square example, we will begin constructing a knot $K^g$ projection on $G^g$ by placing a regular $4g$-gon in the center of $G^g$ with its crossings alternating. Then, we will extend the two strands at each crossing to the corresponding edge of $G^g$ to form $4g$ arms.

The argument that $K^g$ is hyperbolic exactly follows the argument that $K$ is hyperbolic. All complementary regions are open disks, and the central $4g$-gon region and bigon regions are all only adjacent to the outer region. Hence, we follow the same arguments which tell us $K^g$ is cellular alternating and obviously prime. Hence, the knot of which $K^g$ is a projection is hyperbolic.

Next, let $S^g$ be the surface formed by the interior of the central $4g$-gon and each of the bigon arms. Now, we will show $S^g$ is totally geodesic. $S^g$ possesses $4g$-fold rotational symmetry around the center point of the central $4g$-gon. It also possesses reflective symmetry about each of the lines that run through opposite corners of the central $4g$-gon as well as each of the lines that bisect opposite edges of the central $4g$-gon. Additionally, $S^g$ possesses reflective symmetry about the lines corresponding to each edge of $G^g$. Let $L_1$ be one of the lines that runs through two corners of the central square, let $L_2$ be the first edge-bisecting line reached when moving clockwise from $L_1$, and let $L_3$ be the line corresponding to the edge of $G^g$ that $L_1$ intersects. Then, quotienting $S^g$ by the reflections over $L_1, L_2,$ and $L_3$ gives an orbifold $O^g$ which is topologically equivalent to $O$ from before, except it has corner reflectors with orders $4g, 2, \infty,$ and $\infty$. We can once again reduce the two degree $\infty$ corner reflectors to a single point along the edge that connects them. This gives an orbifold $O'^{g}$ with three edges and corner reflectors with orders $4g, 2, \infty$. Since $\frac{1}{4g} + \frac{1}{2} + \frac{1}{\infty} \leq \frac{1}{4 \cdot 2} + \frac{1}{2} = \frac{5}{8} < 1$, $O'^{g}$ is a rigid orbifold. Therefore, $S^g$ is totally geodesic and $K^g$ is a TGS knot.

Given $K^g$ embedded in a thickened genus $g$ surface, we can extend $K^g$ to an infinite family of TGS knots. Take $K^g$ and add $n$ half-twists to each arm of $K^g$, ensuring that the resulting knot is still alternating. Call this knot $K^g_n$. Note, the knot of which $K^g_n$ is a projection is hyperbolic by the same argument as before. Adding twists simply increases the amount of bigons in each arms, but does not change anything else relevant to the hyperbolicity argument. Let $S^g_n$ be the surface consisting of the interior of the central $4g$-gon and each of the twisted arms. $S^g_n$ possesses the same reflective symmetries as $S^g$, and taking the quotient of $S^g_n$ by these lines of reflection gives a rigid orbifold isotopic to $O^g$. Therefore, $S^g_n$ is totally geodesic and $K^g_n$ forms an infinite family of TGS knots.

\subsection{A Second Infinite Family in Each Thickened Surface}

Now, we will show there exists a second infinite family of TGS knots in each thickened surface of arbitrary genus. This second family will be similar yet still distinct from the first family detailed in the previous subsection.

\begin{figure}[htbp]
\begin{center}
\includegraphics[width=4cm]{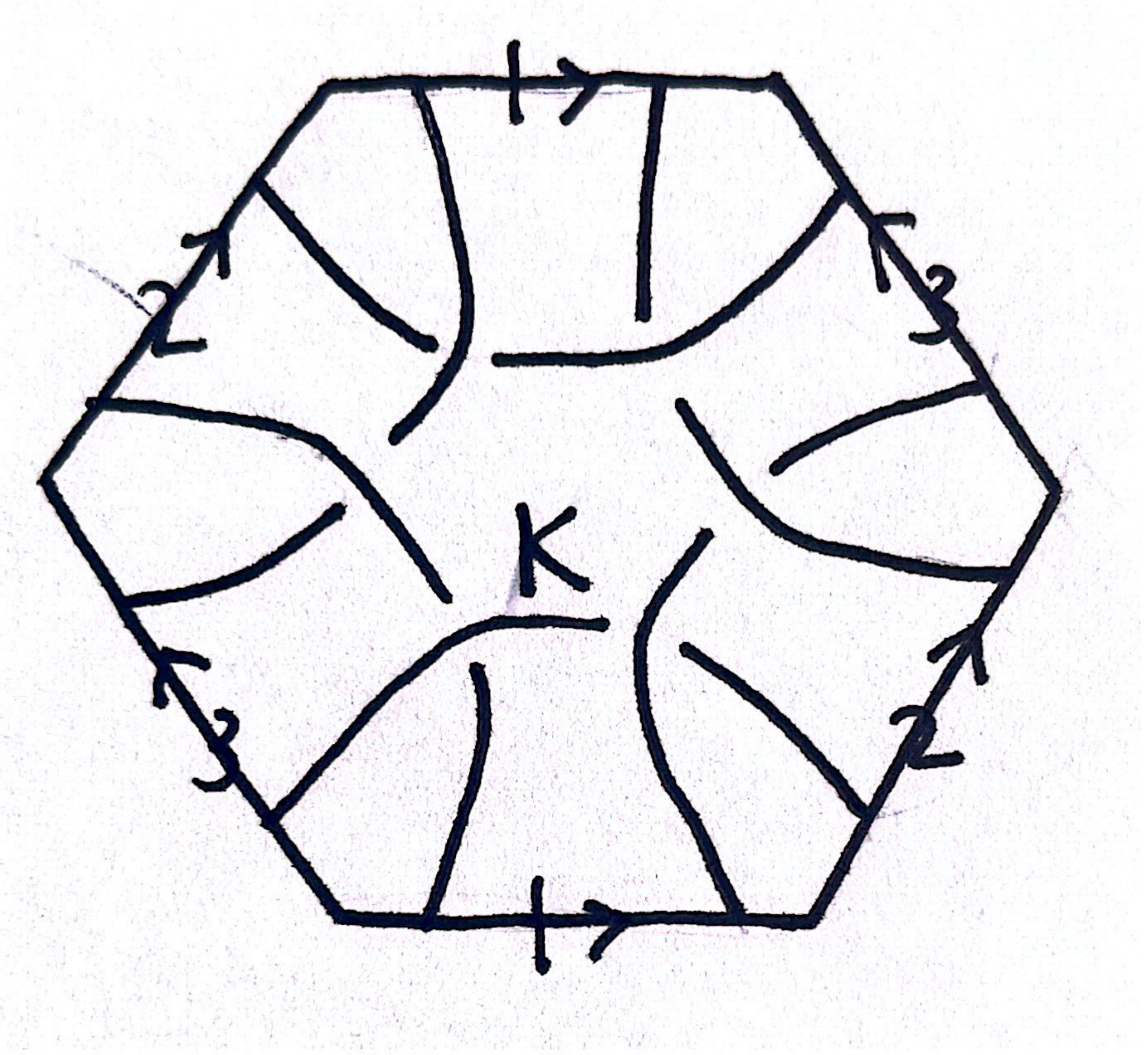}

\caption{}
\label{ThTorHex}
\end{center}
\end{figure}

We will again begin by looking at the thickened torus. Let $\mathcal{G}$ be a regular hexagon gluing diagram with opposite edges identified. Glued together, $\mathcal{G}$ represents a torus $F$. We will consider the thickened surface $F \times [-1, 1]$ with $\mathcal{G} = F \times \{0\}$. We will construct the knot projection $\mathcal{K}$ onto $\mathcal{G}$ shown in Figure \ref{ThTorHex} by placing a regular hexagon in the center and then extending the two strands that meet at each crossing to have one meet the edge on either side of the corresponding vertex of $\mathcal{G}$. Crossings will then be assigned to ensure $\mathcal{K}$ is alternating. 

\begin{figure}[htbp]
\begin{center}
\includegraphics[width=4cm]{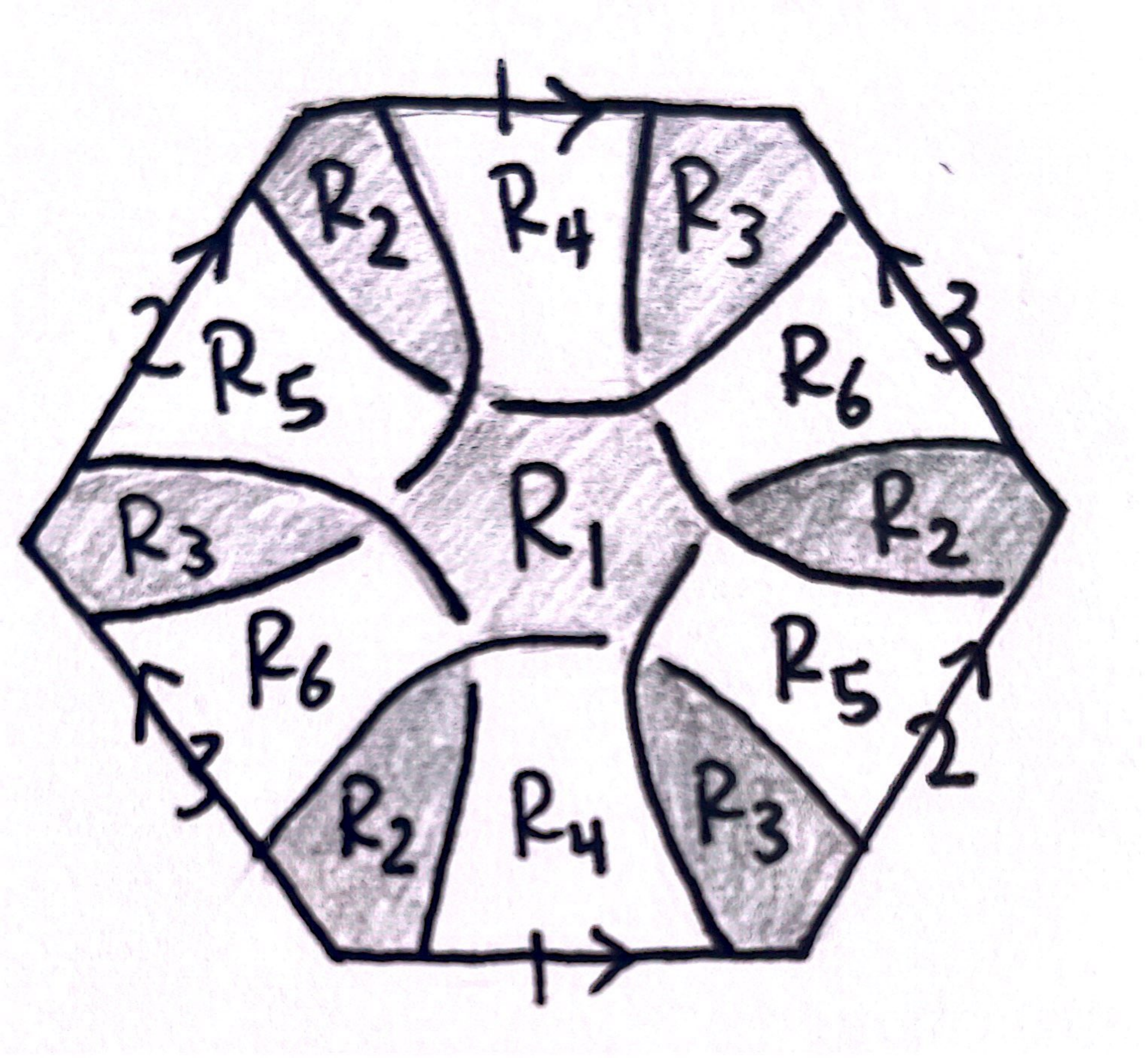}

\caption{}
\label{ThTorHexRegs}
\end{center}
\end{figure}

The complement of $\mathcal{K}$ consists of six regions, each of which has an interior that is a topological open disk. Figure \ref{ThTorHexRegs} depicts these six regions. $\mathcal{R}_1$ is a hexagon, $\mathcal{R}_2$ and $\mathcal{R}_3$ are both triangles, and each of $\mathcal{R}_4, \mathcal{R}_5,$ and $\mathcal{R}_6$ are quadrilaterals. Therefore, $\mathcal{K}$ is cellular alternating.

We can follow the same argument as before to show $\mathcal{K}$ is obviously prime. The only difference in this case is that there are now three outer regions. However, this does not change the argument since $R_1, R_2,$ and $R_3$ are all still only adjacent to the three outer regions. Therefore, we can replace any instance of ``the outer region" or ``$R_4$" in the hyperbolicity argument from Section \ref{th1} with ``an outer region" or ``one of $R_4, R_5, R_6$" and the argument will hold. Hence, the knot of which $\mathcal{K}$ is a projection is hyperbolic in $F \times [-1, 1]$.

\begin{figure}[htbp]
\begin{center}
\includegraphics[width=9cm]{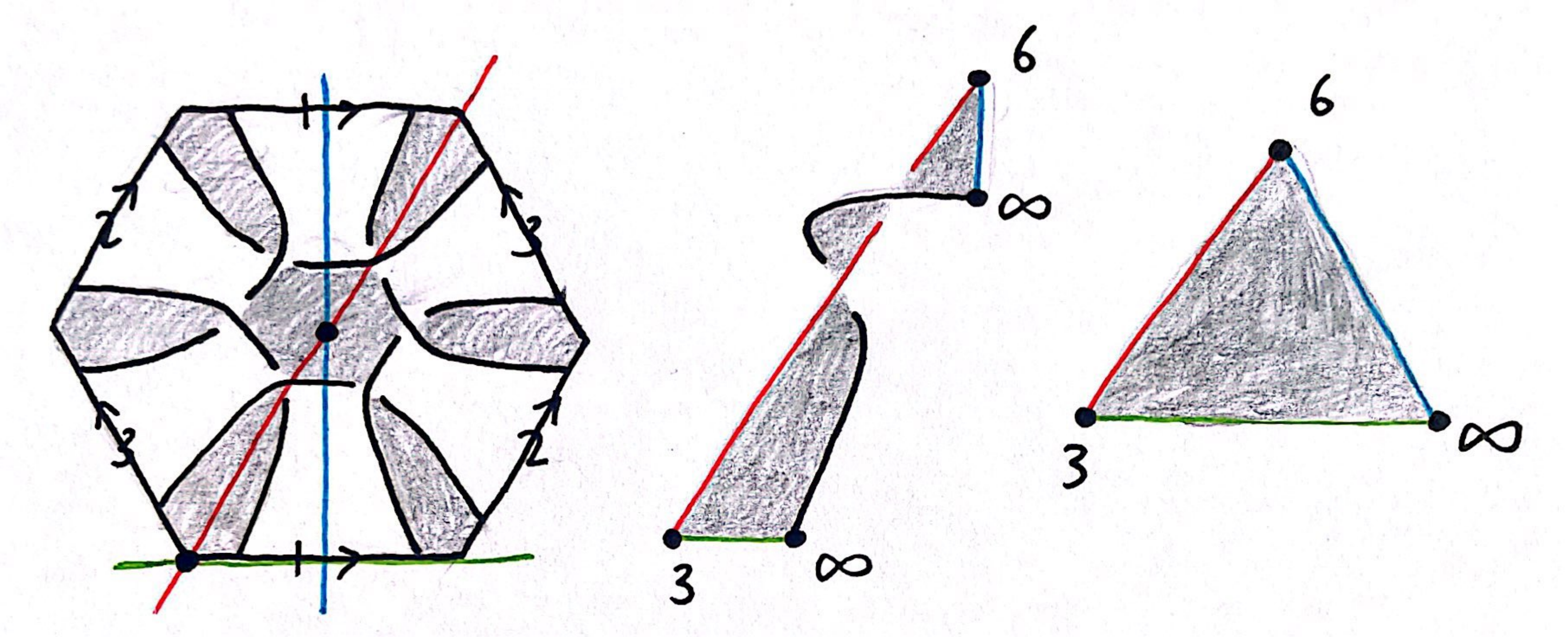}

\caption{}
\label{ThTorHexOrbs}
\end{center}
\end{figure}

We will now show that $\mathcal{K}$ contains a totally geodesic spanning surface by reducing the surface to a rigid orbifold. Let $\mathcal{S} = \mathcal{R}_1 \cup \mathcal{R}_2 \cup \mathcal{R}_3$ be the shaded surface in Figure \ref{ThTorHexRegs}. $\mathcal{S}$ has 6-fold rotational symmetry about the center point of the central hexagon and 3-fold symmetry about both vertices of $\mathcal{G}$. Note, every other vertex of this gluing diagram are identified by the edge identifications, resulting in 2 distinct vertices. Additionally, $\mathcal{S}$ also possesses many reflective symmetries, three of which are shown in the left diagram of Figure \ref{ThTorHexOrbs}. $L_1$, drawn in red, passes through opposite corners of the central hexagon, $L_2$, drawn in blue, passes through the midpoints of opposite edges of the central hexagon, and $L_3$, drawn in green, lays on top of an edge of $\mathcal{G}$. Then, we can take the quotient of $\mathcal{S}$ by these lines of symmetry to get an orbifold $\mathcal{O}$, which is shown in the center diagram of Figure \ref{ThTorHexOrbs}. As before, the two degree $\infty$ corner reflectors in $\mathcal{O}$ can be reduced to a single corner reflector along the edge of $\mathcal{O}$ corresponding to a segment of $\mathcal{K}$. This gives an equivalent orbifold $\mathcal{O}'$ with corner reflectors of orders $6, 3,$ and $\infty$. $\mathcal{O}'$ is shown in the right diagram of \ref{ThTorHexOrbs}. Since $\frac{1}{6} + \frac{1}{3} + \frac{1}{\infty} = \frac{1}{2} < 1$, it follows by \cite{orbifolds} that $\mathcal{O'}$ is a rigid orbifold. Therefore, $\mathcal{S}$ is totally geodesic and $\mathcal{K}$ is a TGS knot.

Next, we will extend $\mathcal{K}$ to create an infinite family of TGS knots. Begin with $\mathcal{K}$ and add $n$ full-twists to each arm so that the resulting knot is still alternating, and let this new knot be $\mathcal{K}_n$. Note, we cannot add half-twists to each arm of $\mathcal{K}$ as that would instead result in a link rather than a knot. 

As we discussed in Section \ref{th1}, adding twists does not alter the arguments for the projection being cellular alternating and obviously prime. Hence, we can use the same argument as $\mathcal{K}$ and we get that the knot of which $\mathcal{K}_n$ is a projection is hyperbolic in $F \times [-1,1]$.

\begin{figure}[htbp]
\begin{center}
\includegraphics[width=4cm]{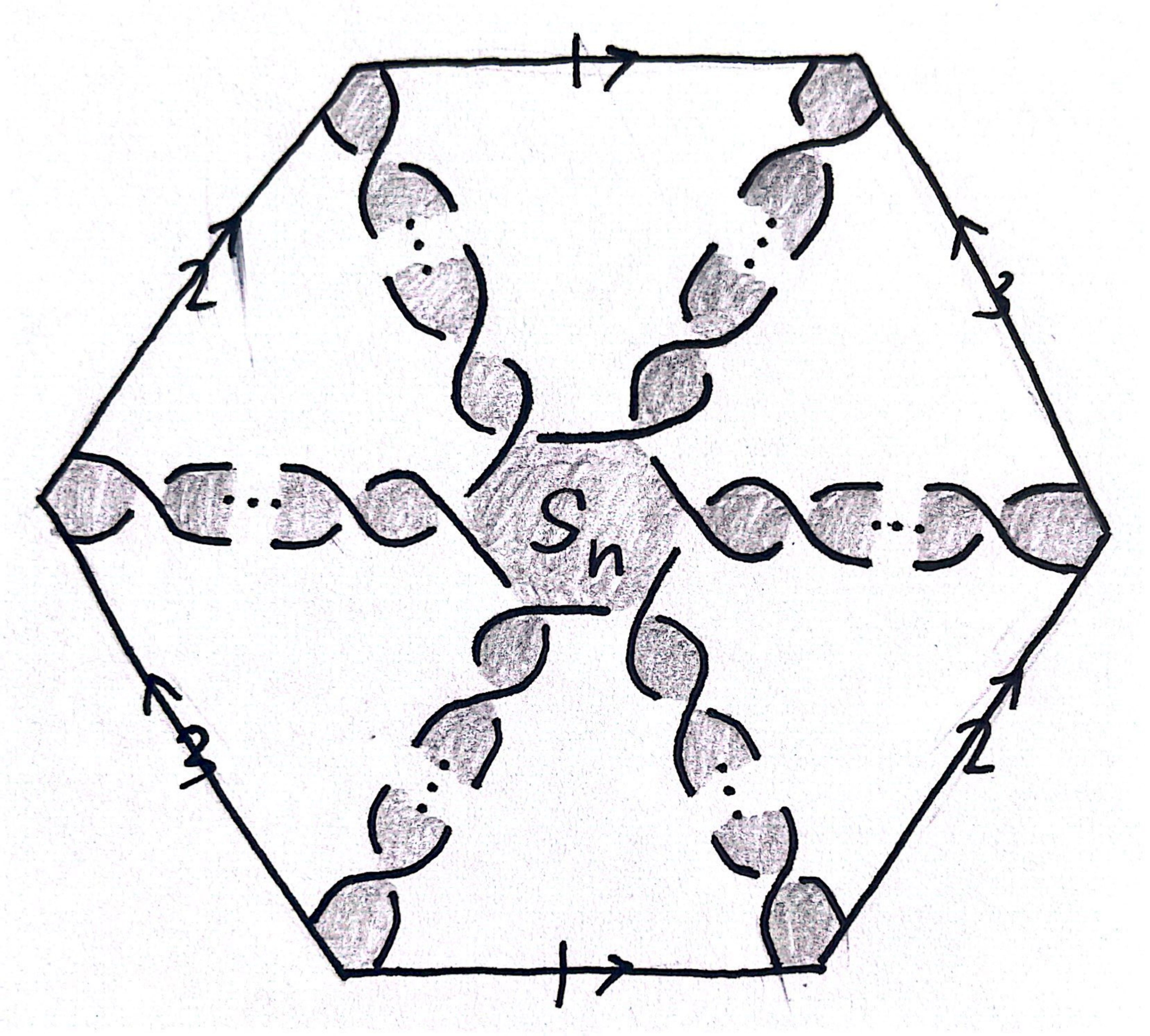}

\caption{}
\label{ThTorHexKn}
\end{center}
\end{figure}

Let $\mathcal{S}_n$ be the surface formed by the interior of the central hexagon and each of the twisted arms. Figure \ref{ThTorHexKn} depicts $\mathcal{S}_n$. 

\begin{figure}[htbp]
\begin{center}
\includegraphics[width=9cm]{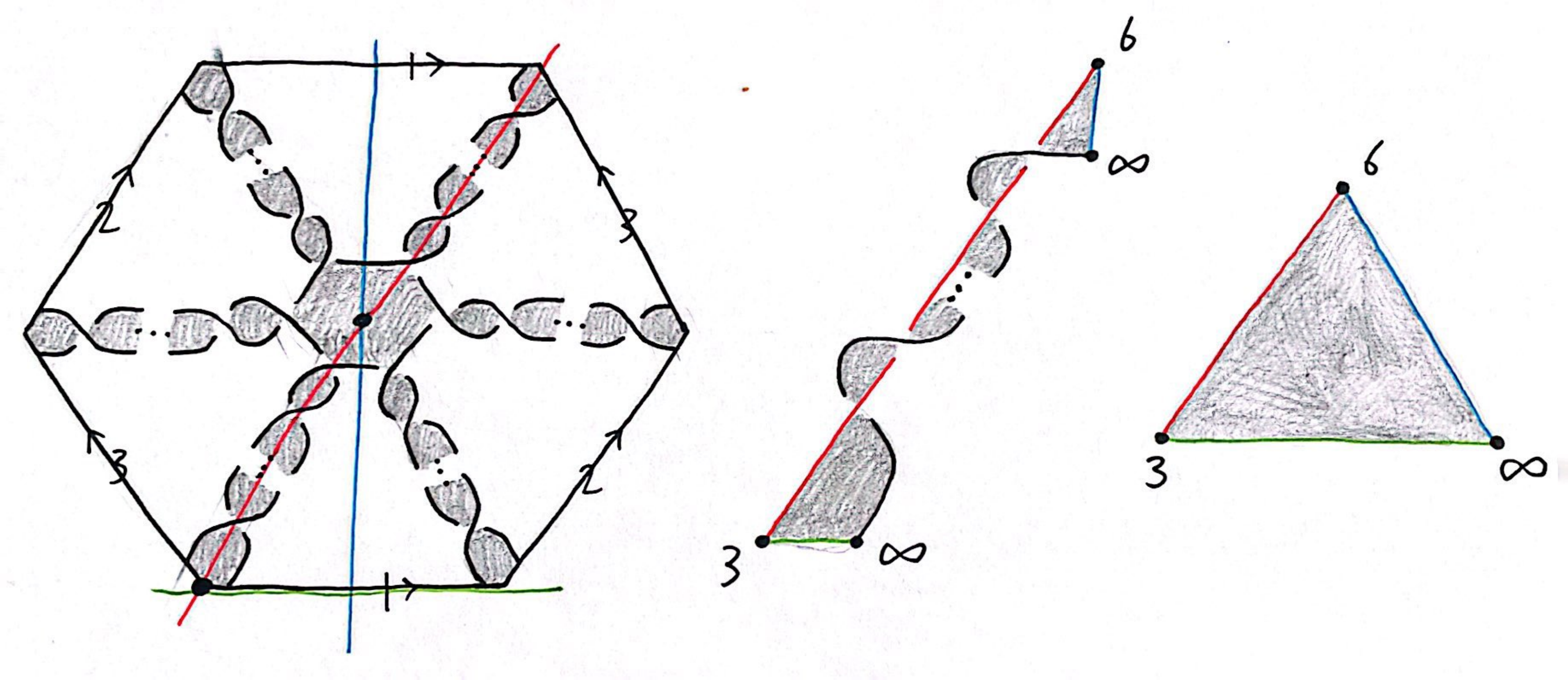}

\caption{}
\label{ThTorHexKnOrbs}
\end{center}
\end{figure}

We will now show that $\mathcal{S}_n$ is totally geodesic. $\mathcal{S}_n$ has reflective symmetry about $L_1, L_2,$ and $L_3$ just like $\mathcal{S}$. These lines of symmetry are illustrated in the left diagram of Figure \ref{ThTorHexKnOrbs}. Quotienting $\mathcal{S}_n$ by these symmetries gives an orbifold $\mathcal{O}_n$, shown in the center diagram of Figure \ref{ThTorHexKnOrbs}, which is isotopic to $\mathcal{O}$ from before. Therefore, $\mathcal{O}_n$ is a rigid orbifold. Hence, $\mathcal{S}_n$ is totally geodesic and all of the $\mathcal{K}_n$ form an infinite family of TGS knots in the thickened torus.

In the last part of this section, we will now show that this infinite family of TGS knots can be extended to thickened surfaces of any genus $g$ with $g \geq 2$. Let $\mathcal{G}^g$ be a $(4g + 2)$-gon gluing diagram  with opposite edges identified. Then, $\mathcal{G}^g$ represents a genus $g$ surface $\mathcal{F}^g$. Consider $\mathcal{F}^g \times [-1,1]$ with $\mathcal{G}^g = \mathcal{F}^g \times \{0\}$. To construct a knot $\mathcal{K}^g$ embedded in $\mathcal{G}^g$, place a regular $(4g+2)$-gon in the center of $\mathcal{G}^g$ and extend the two strands of each crossing to have one meet the two edges one either side of the corresponding vertex of $\mathcal{G}^g$.

The proof that the knot of which $\mathcal{K}^g$ is a projection is hyperbolic is the same as that for $\mathcal{K}$ using Theorem \ref{small2018}. Let $\mathcal{S}^g$ be the surface formed by the interior of the central $(4g+2)$-gon and the interiors of each of the arms. Therefore, $\mathcal{S}^g$ consists of one $(4g+2)$-gon and two $(2g+1)$-gons. $\mathcal{S}^g$ possesses $(4g+2)$ fold symmetry about the center of the central $(4g+2)$-gon and $(2g+1)$-fold symmetry about each vertex of $\mathcal{G}^g$. In addition, $\mathcal{S}^g$ has reflective symmetry about each line that connects opposite vertices or midpoints of opposite edges of $\mathcal{G}^g$ as well as each line corresponding to an edge of $\mathcal{G}^g$. Let $L_1$ be a line of reflective symmetry which passes through two vertices of $\mathcal{G}^g$, let $L_2$ be the the next line of symmetry which connects edge midpoint that is counter-clockwise from $L_1$, and let $L_3$ be the line of symmetry corresponding to the edge of $\mathcal{G}^g$ which $L_2$ intersects. Taking the quotient of $\mathcal{S}^g$ by these three symmetries gives an orbifold $\mathcal{O}^g$ equivalent to $\mathcal{O}$, except with corner reflectors of orders $4g+2, 2g+1, \infty,$ and $\infty$. The two degree $\infty$ corner reflectors can be reduced again to a single point, which gives a triangular orbifold. Note, $\frac{1}{4g+2} + \frac{1}{2g+1} + \frac{1}{\infty} \leq \frac{1}{6} + \frac{1}{3} = \frac{1}{2} < 1$. Therefore, $\mathcal{O}^g$ is rigid. It follows that $\mathcal{S}^g$ is totally geodesic and $\mathcal{K}^g$ is a TGS knot.

Lastly, $\mathcal{K}^g$ can be extended to an infinite family of TGS knots by adding $n$ full-twists to each arm. Let this new knot in $F^g \times [-1,1]$ be $\mathcal{K}^g_n$. Again, the knot of which $\mathcal{K}^g_n$ is a projection is hyperbolic by the same argument as before. Let $\mathcal{S}^g_n$ be the surface formed by the interior of the central $(4g+2)$-gon and the interiors of each arm of $\mathcal{K}^g_n$. $\mathcal{S}^g_n$ has the same symmetries as $\mathcal{K}^g$, and quotienting $\mathcal{S}^g_n$ by these symmetries gives a rigid orbifold which is isotopic to $\mathcal{O}^g$. Hence, $\mathcal{S}^g_n$ is totally geodesic and all of the $\mathcal{K}^g_n$ form an infinite family of TGS links.

Every TGS knot discussed in this section correspond to a virtual knot, just like our first TGS knot $K$. Therefore, every totally geodesic spanning surface mentioned in this section corresponds to a totally geodesic spanning surface of a virtual knot.

\section{TGS Links in the Sphere Cross the Circle}\label{SXC}

\begin{figure}[htbp]
\begin{center}
\includegraphics[width=6.5cm]{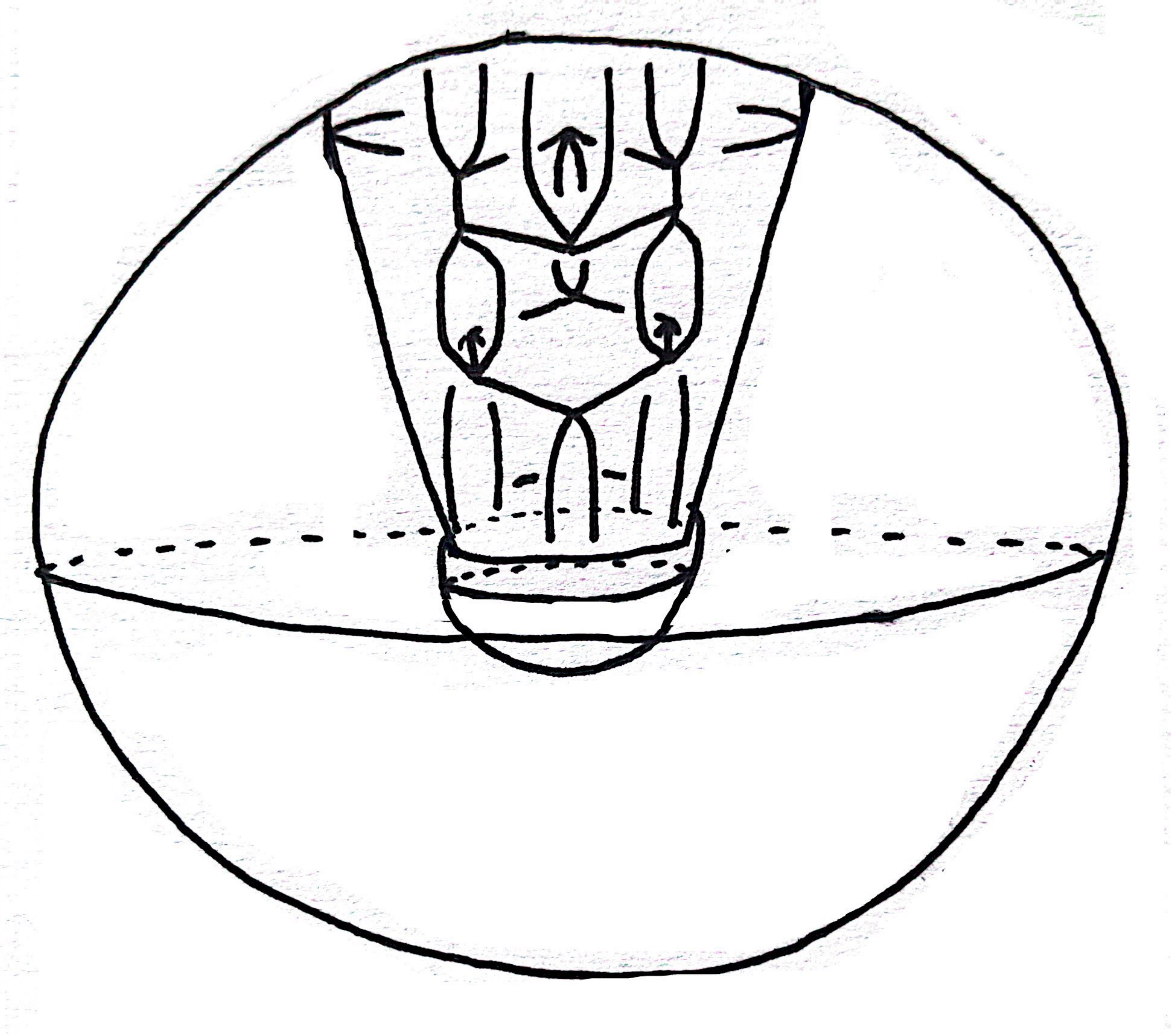}

\caption{}
\label{SphXCir}
\end{center}
\end{figure}

In this section, we will look at infinite families of TGS links embedded in $\mathcal{S} = S^2 \times S^1$. We will begin by depicting $\mathcal{S}$ as $S^2 \times [0,1]$ such that $S^2 \times \{0\}$, the outer sphere, is identified with $S^2 \times \{1\}$, the inner sphere. We will now construct a six-component link $L$ which is the simplest in a family of links we will call ``layer cake" links. Begin by embedding two hexagons of link strands at $S^2 \times \{\frac{1}{4}\}$, call this hexagon $P$, and $S^2 \times \{\frac{3}{4}\}$, call this hexagon $Q$. Label the vertices of the $P$ as $v_1, v_2, v_3, v_4, v_5, v_6$ and label the vertices of $Q$ as $w_1, w_2, w_3, w_4, w_5, w_6$ so that $v_i$ lies directly above $w_i$. Then, connect $v_1$ to $w_1$, $v_3$ to $w_3$, and $v_5$ to $w_5$ using bigon ``arms" that run through $S^2 \times [\frac{1}{4}, \frac{3}{4}]$. Next, connect $v_2$ to $w_2$, $v_4$ to $w_4$, and $v_6$ to $w_6$ using bigons arms that run through $[0,\frac{1}{4}] \cup [\frac{3}{4}, 1]$. A depiction of $L$ can be seen in Figure \ref{SphXCir}.

\begin{figure}[htbp]
\begin{center}
\includegraphics[width=5.5cm]{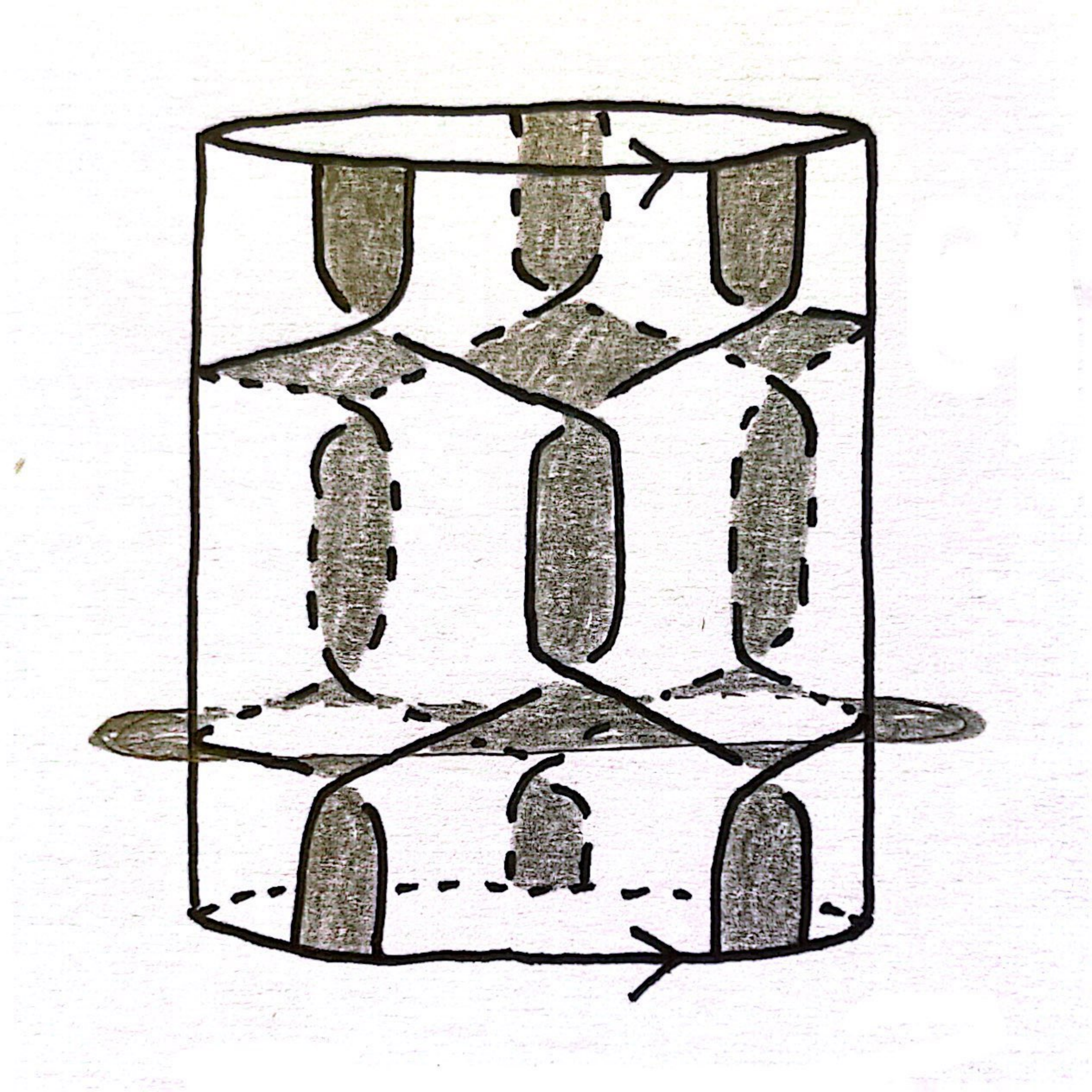}

\caption{}
\label{SXCTor}
\end{center}
\end{figure}

We will then project $L$ out to the depicted cylinder which surrounds it. This cylinder forms a torus, whose boundary we will call $\mathcal{T}$, since the identification of the inner and outer spheres identifies the top to the bottom of this cylinder. Note, $\mathcal{S}$ can also be constructed from gluing two tori $T_1, T_2$ together along their boundaries. Then, $\mathcal{T}$ is these glued boundaries, and the inside and outside of $\mathcal{T}$ represent $T_1$ and $T_2$. Now, we can assign crossings to each intersection so that $L$ is alternating. This projection can be seen in Figure \ref{SXCTor}.

Any curve which bounds a compressing disk of $\mathcal{T}$ must necessarily wrap once meridinally around $\mathcal{T}$. Therefore, we will find a lower bound on the representativity of $L$ by counting the minimum number of intersections any such curve must have with $L$. Beginning at some point on $\mathcal{T}$ not on $L$, we will travel clockwise until we reach an arm of $L$. To get past this arm, our first option is to travel directly through the arm. In this case, the curve will intersect the two strands which form the arm. Therefore, the curve will have intersected $L$ two times. The other option would be to go around the arm by crossing over a strand of a layer polygon and entering a layer adjacent to the one in which we began. This will cause the curve to intersect $L$ one time. Additionally, in order eventually close the curve, we must travel through some strand of that layer polygon again to re-enter the layer in which we began. Therefore, this option also causes the curve to intersect $L$ two times. Since each layer of $L$ contains three arms, any curve which bounds a compressing disk must intersect $L$ at least six times. Hence, the representativity of $L$ is at least 6. It follows by Theorem \ref{HowiePurcell} that $L$ is hyperbolic in $\mathcal{S}$ .

\begin{figure}[htbp]
\begin{center}
\includegraphics[width=13cm]{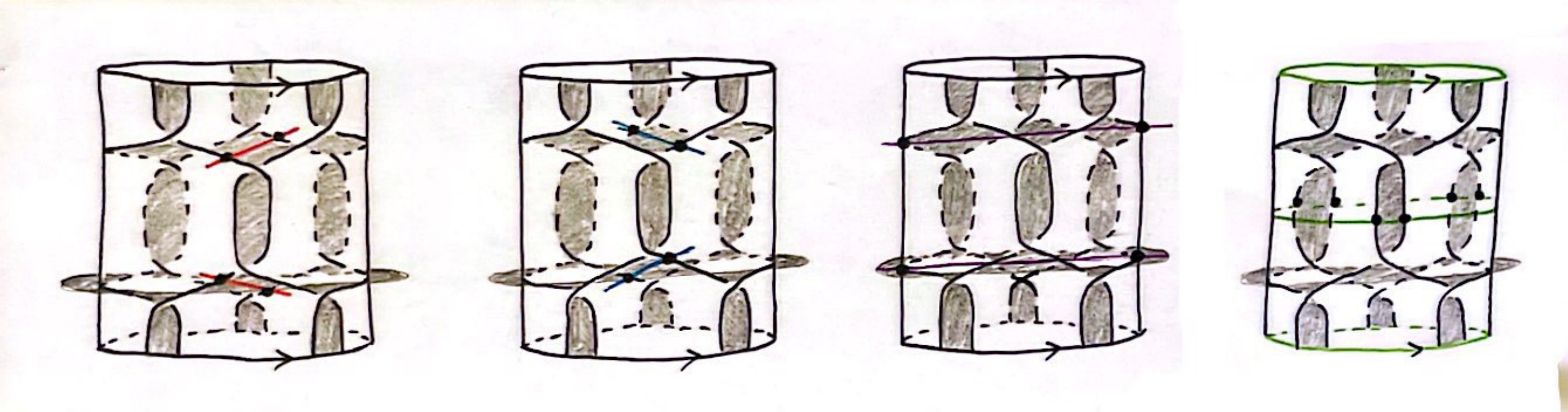}

\caption{}
\label{SXCTorSym}
\end{center}
\end{figure}

Now, consider the spanning surface $S$ of $L$ composed of the interior of $P$ in $T_1$, the interior of $Q$ in $T_2$, and the interiors of each bigon arm along $\mathcal{T}$. This can be seen in Figure \ref{SXCTor}. We will show that $S$ is totally geodesic by quotienting $S$ down to a rigid orbifold. First, $S$ possesses rotational symmetry by $180^\circ$ about the four pairs of lines and meridian circles shown in Figure \ref{SXCTorSym}. The two red lines $L_{11}, L_{12}$ in the first diagram of Figure \ref{SXCTorSym} lie in $S^2 \times \{\frac{1}{4}\}$ and $S^2 \times \{\frac{3}{4}\}$ and they are perpendicular bisectors of $\overline{v_1v_2}, \overline{v_4v_5} \in P$ and $\overline{w_1w_2}, \overline{w_4w_5} \in Q$, respectively. The two blue lines $L_{21}, L_{22}$ in the second diagram of Figure \ref{SXCTorSym} lie in $S^2 \times \{\frac{1}{4}\}$ and $S^2 \times \{\frac{3}{4}\}$ and they are perpendicular bisectors of $\overline{v_2v_3}, \overline{v_5v_6} \in P$ and $\overline{w_2w_3}, \overline{w_5w_6} \in Q$, respectively. The two purple lines $L_{31}, L_{32}$ in the third diagram of Figure \ref{SXCTorSym} lie in $S^2 \times \{\frac{1}{4}\}$ and $S^2 \times \{\frac{3}{4}\}$ and they are perpendicular bisectors of $\overline{v_3v_4}, \overline{v_1v_6} \in P$ and $\overline{w_3w_4}, \overline{w_1w_6} \in Q$, respectively. Lastly, the two green circles $C_{41}, C_{42}$ in the fourth diagram of Figure \ref{SXCTorSym} lie in $S^2 \times \{\frac{1}{2}\}$ and $S^2 \times \{0\} = S^2 \times \{1\}$ and they are perpendicular bisectors of the three arms which run from $v_1$ to $w_1$, $v_3$ to $w_3$, and $v_5$ to $w_5$ and the three arms which run from $v_2$ to $w_2$, $v_4$ to $w_4$, and $v_6$ to $w_6$, respectively. The red, blue, and purple symmetries send the link segments and associated surface sections in $S^2 \times [\frac{1}{4}, \frac{3}{4}]$ to $S \times [0, \frac{1}{4}] \cup [\frac{3}{4}, 1]$ and vice versa. The green symmetry swaps these two regions and it also sends $T_1$ to $T_2$.

\begin{figure}[htbp]
\begin{center}
\includegraphics[width=7cm]{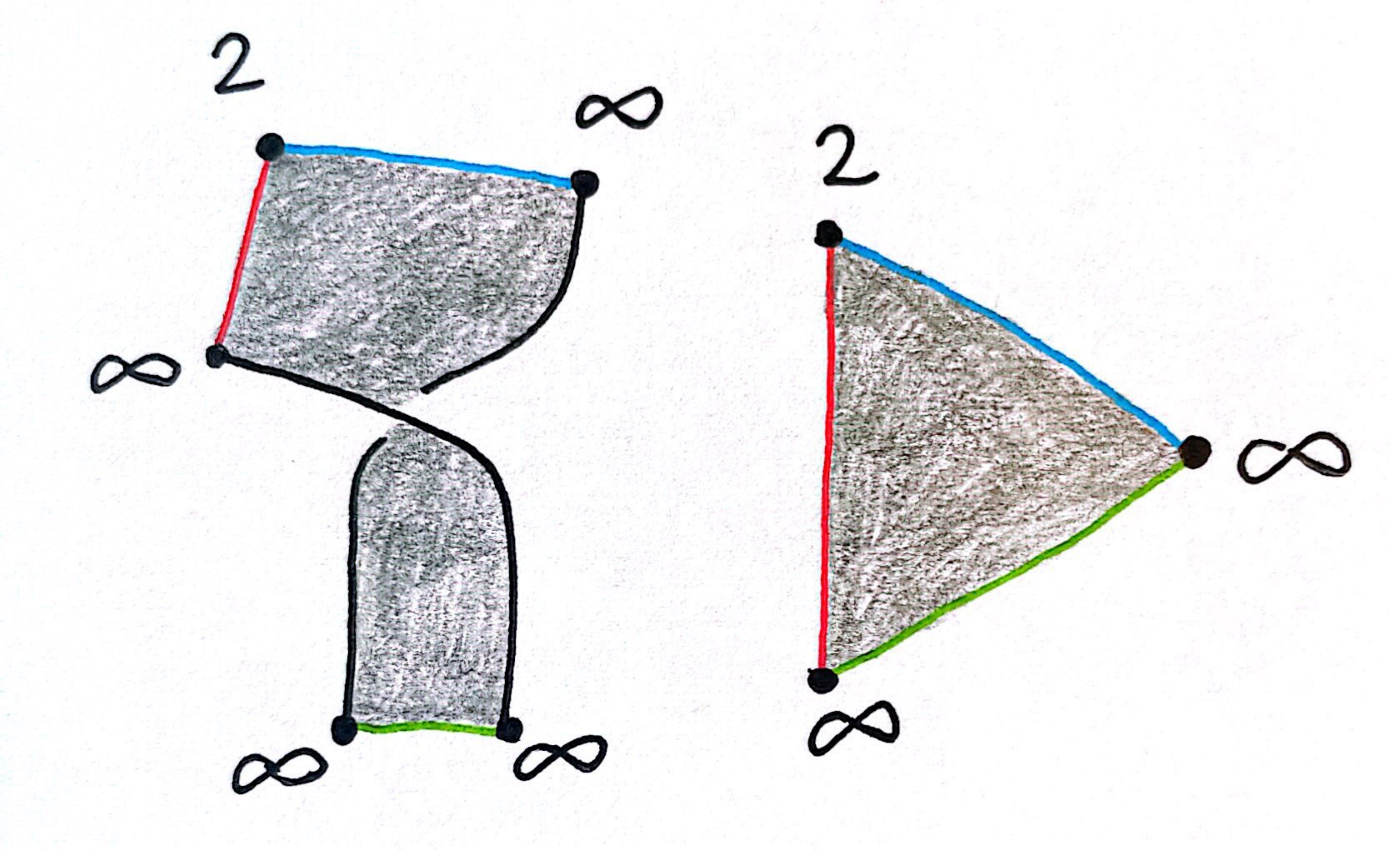}

\caption{}
\label{SXCTorOrb}
\end{center}
\end{figure}

Taking the quotient of $S$ by the red, blue, and green symmetries gives an orbifold $O$ as in the left diagram of Figure \ref{SXCTorOrb}. $O$ possesses five corner reflectors with orders $2, \infty, \infty, \infty,$ and $\infty$. However, as we have seen before, two pairs of order $\infty$ corner reflectors can be reduced to a point along the sides of $O$ which are segments of $L$. The yields an equivalent orbifold $O'$, seen in the right diagram of \ref{SXCTorOrb}. Since $O'$ has three corner reflectors of orders $2, \infty,$ and $\infty$, and $\frac{1}{2} + \frac{1}{\infty} + \frac{1}{\infty} = \frac{1}{2} < 1$, $O'$ is a rigid orbifold by the result from \cite{orbifolds}. Hence, $S$ is totally geodesic and $L$ is a TGS link in $\mathcal{S}$. 

$L$ has multiple ways it can be extended to form an infinite family of TGS links. This includes adding more hexagons and connecting arms to create more layers, adding sides to $P$ and $Q$ and then adding the appropriate connecting arms, and adding twists to each arm.

First, we will consider adding more layers. Begin with $L$ and add more pairs of layers to the ``layer cake" so that there are a total of $2^\ell$ layers. That is, add $2^\ell$ hexagons and connect them to each other and the $L$, making sure to adhere to the pattern of connecting every other arm to the hexagon above and the hexagon below. Call this new link $L'_\ell$. Since we still have that every layer of $L'_\ell$ has three arms, the representativity of $L'_\ell$ is at least 6. Therefore, $L'_\ell$ is hyperbolic. Let $S'_\ell$ be the surface formed by the interiors of each arm of $L'_\ell$ on $mathcal{T}$ and the interior of each hexagon in either $T_1$ or $T_2$, alternating between layers. Now, we will show that $S'_\ell$ can be quotiented down so that it becomes $S$. Let $C_1$ and $C_2$ be circles which lie in $S^2 \times \{\frac{1}{2}\}$ and $S^2 \times \{1\}$ and bisect three arms of $S'_\ell$. Then, $S'_\ell$ has the same rotational symmetry about $C_1$ and $C_2$ as we saw before. Quotienting $S'_\ell$ by these symmetries yields a layer cake link and corresponding surface that has $2^{\ell-1}$ layers. We can repeat this process until we are left with a surface which has only 2 layers. That is, until we are left with $S$. Therefore, $S'_\ell$ can be further reduced to form a rigid orbifold, and thus it is totally geodesic. It also follows that $L'_\ell$ is a TGS link.

Next, we will look at adding sides to the polygon in each layer. Once again, begin with $L$. This time, add edges to each of $P$ and $Q$ to form $P''$ and $Q''$ such that each have $2m$ edges, and add arms which connect the two vertices so that arms between odd-indexed vertices run through $S^2 \times [\frac{1}{4}, \frac{3}{4}]$ and arms between even-indexed vertices run through $S \times [0, \frac{1}{4}] \cup [\frac{3}{4}, 1]$. Call this new link $L''_m$. Increasing the number of arms will increase the lower bound on the representativity of $L''_m$, which means it remains greater than 4. Therefore, $L''_m$ is hyperbolic. Let $S''_m$ be the surface made from the interiors of $P''$ in $T_1$, $Q''$ in $T_2$, and each of the arms on $\mathcal{T}$. Note, $S''_m$ still possesses rotational symmetry about the circles $C_{31}$ and $C_{32}$ from the surface $S$ case. Additionally, $S''_m$ has rotational symmetry about each pair of lines which run through corresponding sets of opposite edges of $P''$ and $Q''$. These lines are analogous to $L_{11}, L_{12}, L_{21},$ and $L_{22}$ from the surface $S$ case. Therefore, $S''_m$ can be quotiented down by two pairs of these lines which run through adjacent edges of $P''$ and $Q''$ in addition to $C_{31}$ and $C_{32}$ in order to get the orbifold $O$ from before. It follows that $S''_m$ is totally geodesic and $L''_m$ is a TGS link.

Lastly, we will consider adding twists to each arm. Take $L$ and add $n$ half-twists to each arm. Consider the case where $n$ is even. Let this new knot be $L'''_n$, and let the surface formed by the interiors of $P$ in $T_1$, $Q$ in $T_2$, and each twisted arm on $\mathcal{T}$ be $S'''_n$. This case closely resembles the case of twisted arms from the Section \ref{thick}. Since each arm remains in the vertical plane it was in prior to twisting and the twisting has no effect on $P$ or $Q$, $S'''_n$ retains the rotational symmetries depicted in Figure \ref{SXCTorSym} from before. Therefore, we can take the quotient of $S'''_n$ by these symmetries, and then use isotopy moves to untwist the half-arm in the resulting orbifold to get $O$ which is rigid. Hence, $S'''_n$ is totally geodesic and $L'''_n$ is a TGS link.

\begin{figure}[htbp]
\begin{center}
\includegraphics[width=3.5cm]{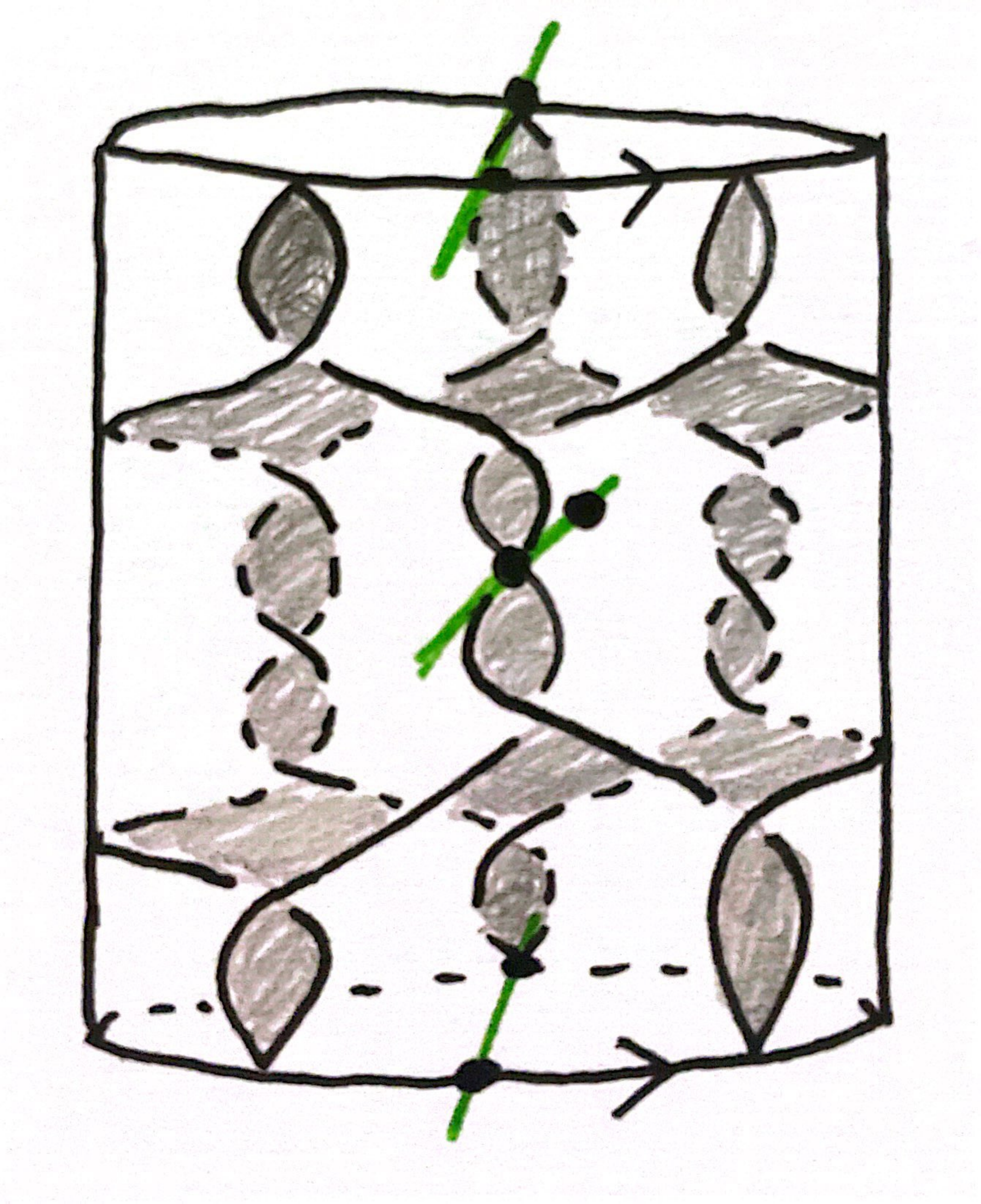}

\caption{}
\label{SXCOdd}
\end{center}
\end{figure}

Next, consider the case where $n$ is odd. Again, let this knot be $L'''_n$. Now, let $S'''_n$ be the surface formed by the interiors $P$ and $Q$ in $T_1$ and each arm on $\mathcal{T}$. Note, since we added an odd number of twists, each arm will now have a crossing at its midpoint. Let $L_{41}$ and $L_{42}$ be lines in the same vertical plane which pierce the arms connecting $v_1$ to $w_1$ and $v_4$ to $w_4$. Then, $S'''_n$ has $180^\circ$ rotational symmetry about this pair of lines, which are depicted in Figure \ref{SXCOdd} for $n=1$. Again, we can take the quotient of $S'''_n$ by these symmetries, and then use isotopy moves to untwist the half-arm in the resulting orbifold to get $O$. Therefore, $S'''_n$ is totally geodesic and $L'''_n$ is a TGS link.

Now, we can combine these three alterations to make a single infinite family of layer cake TGS links. Let $L_{\ell m n}$ be a layer cake link with $2\ell$ layers, each of which has a $(2m)$-gon, and layers are connected by arms with $n$ half-twists. Then, let $S_{\ell m n}$ be the surface formed by the interiors of all connecting arms on $\mathcal{T}$ and the interiors of each $(2m)$-gon, either alternating between these interiors being taken in $T_1$ and $T_2$ or having all interiors being taken in $T_1$ based on the parity of $n$. To show that $S_{\ell m n}$ is totally geodesic, we can follow each of our reduction strategies to remove the added layers, edges, and twists. First, repeatedly quotient $S_{\ell m n}$ by the pair of circles which run through groups of arms until there are only two layers remaining. Then, quotient the resulting surface by two pairs of lines which run through adjacent pairs of edges of the $(2m)$-gons and the pair of circles or lines which bisect arms. Lastly, use isotopy moves to untwist the half-arm in order to get $O$, a rigid orbifold. Hence, $S_{\ell m n}$ is totally geodesic and all links of the form $L_{\ell m n}$ are TGS links in $\mathcal{S}$.

\section{TGS Links in Lens Spaces}

In this section we will look at two separate families of TGS links in lens spaces.

\subsection{Layer Cake Links}

\begin{figure}[htbp]
\begin{center}
\includegraphics[width=10cm]{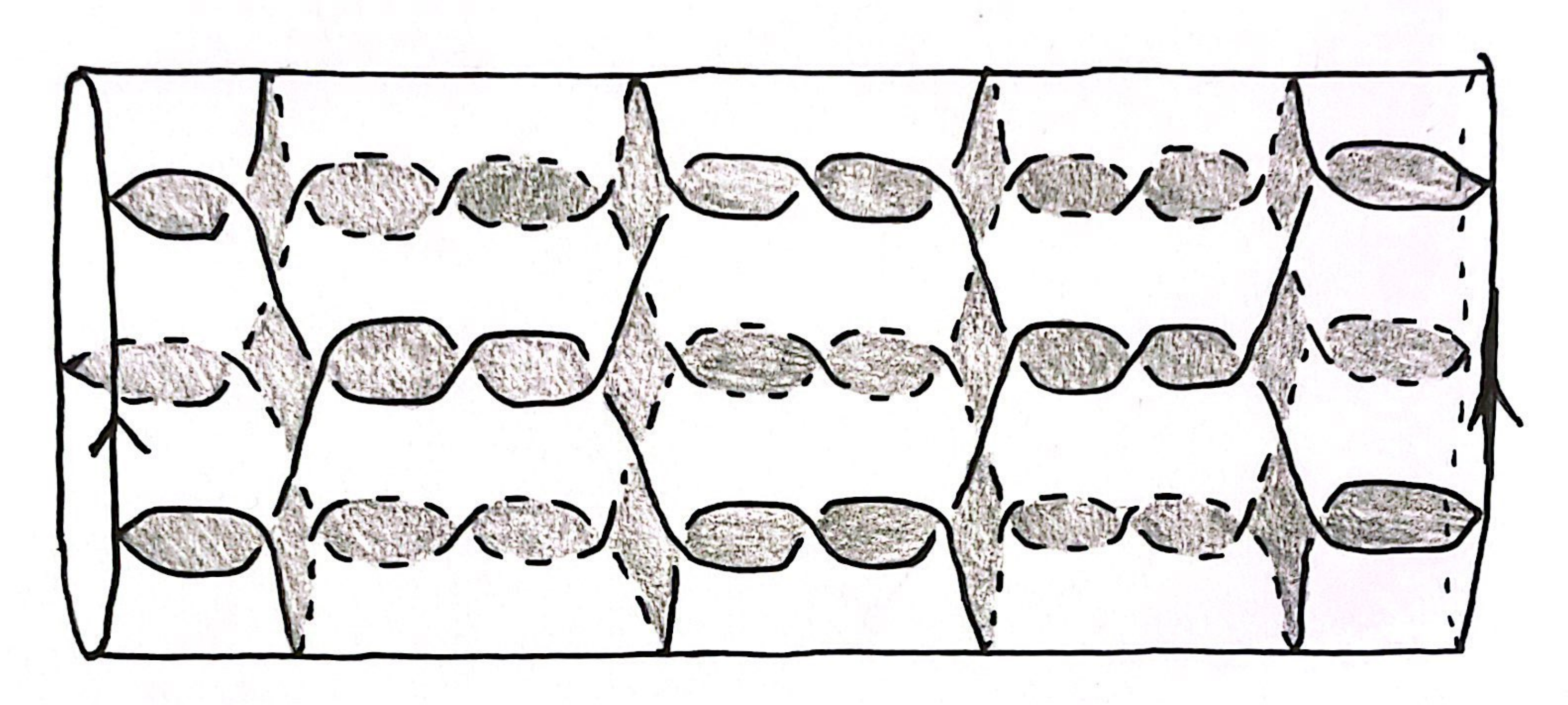}

\caption{}
\label{LensLC}
\end{center}
\end{figure}

Let $L(m,n)$ with $m,n$ relatively prime be a lens space. $L(m,n)$ can be created by gluing two solid tori together along their boundaries such that a meridian curve $c_1$ on the boundary of the first solid torus $T_1$ is sent to an $(n,m)$-curve $c_2$ on the boundary of the second solid torus $T_2$. Note, all meridian curves of $T_1$ parallel to $c_1$ are sent to $(n,m)$-curves parallel to $c_2$ by this gluing. If $n \geq 3$, then onto the boundary of $T_2$ we will project a ``layer cake" link $L$ with $2m$ copies of a $2n$-gon connected by bigon arms. We will refer to these $2n$-gons as ``layer polygons." If $n = 1$ or $n=2$, we will replace the $2n$-gon layer polygons with octagons. We will also add a half-twist to each arm. Lastly, we will assign crossings to each intersection so that $L$ is alternating. Figure \ref{LensLC} depicts $L$ for the lens space $L(2,3)$.

Since $n \geq 3$, each layer of $L$ possesses at least three arms. Therefore, by the same reasoning from Section \ref{SXC}, the representativity of $L$ is at least 6. Then, $L$ is hyperbolic in $L(m,n)$ by Theorem \ref{HowiePurcell}.

\begin{figure}[htbp]
\begin{center}
\includegraphics[width=6cm]{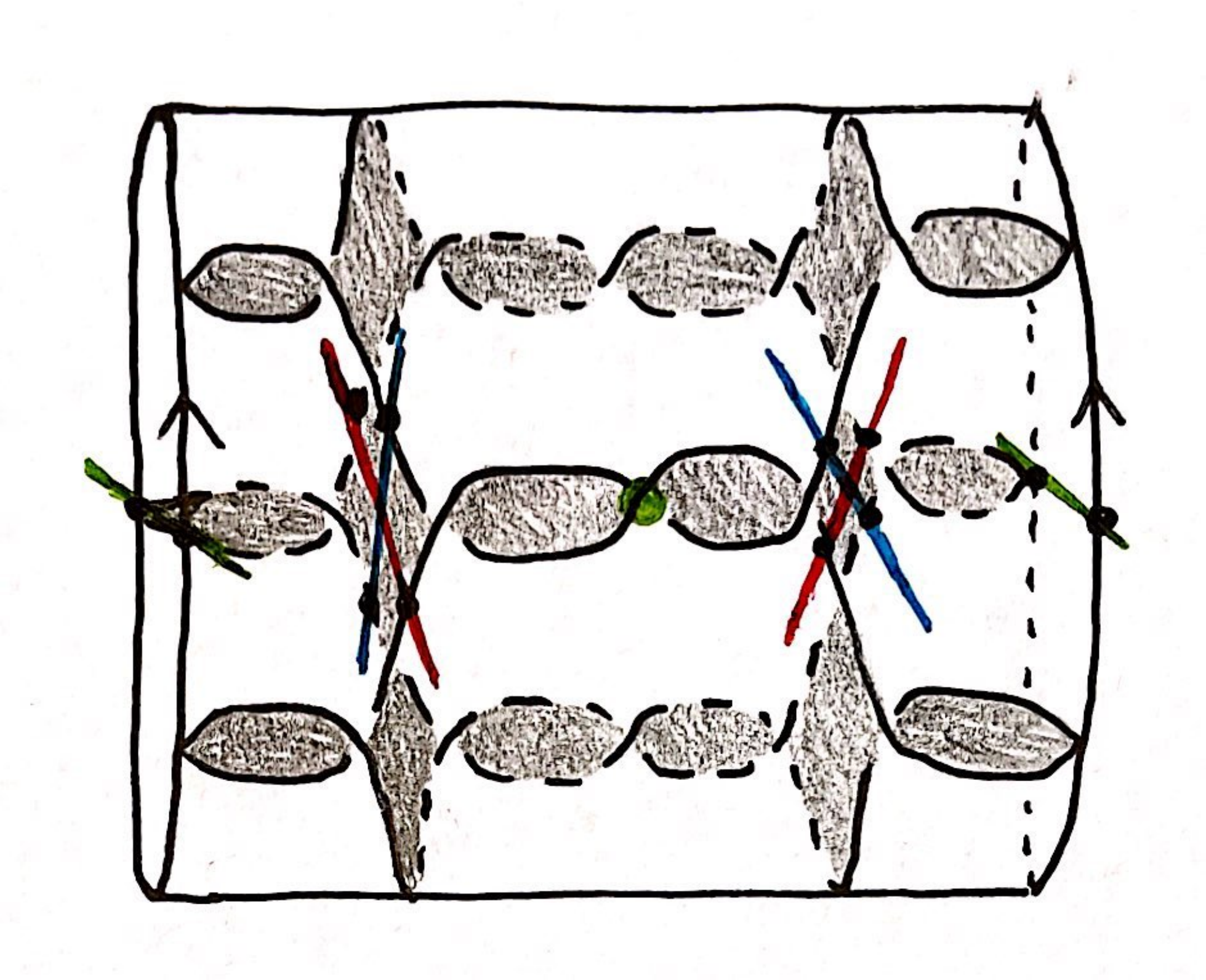}

\caption{}
\label{LensLCSym}
\end{center}
\end{figure}

Next, let $S$ be the surface formed in $T_2$ by the interiors of each layer polygon and each arm. A $\frac{1}{n}$ meridinal and $\frac{1}{m}$ longitudinal rotation of $T_2$ will send $L$ to itself, $S$ to itself, and send any $(n,m)$ curve to itself. Call this symmetry $R$. Therefore, $R$ is a symmetry of $L(m,n)$, which means we can quotient $S$ by this symmetry. Doing so leaves us with a surface $S'$ that consists of a pair of layer polygons connected by bigon arms. Then, $S'$ has the $180^\circ$ rotational symmetries about lines in the layer polygons which bisect opposite edges of the layer polygons as discussed in Section \ref{SXC}. Additionally, $S'$ possesses $180^\circ$ rotational symmetry about lines that pierce the central crossings in opposite arms like in the case in Section \ref{SXC} where an odd number of half twists were added to each arm. Then, we can take two rotational symmetries whose axes of rotation bisect adjacent edges of some layer polygon along with a symmetry that pierces the arm formed by the link segments which make those edges, and we can quotient $S'$ by these three symmetries. One such trio of symmetry axes is shown in Figure \ref{LensLCSym} drawn in red, blue, and green. Note, we are peering down the green axis in this figure, so it appears as a single point. Taking this quotient will give the same orbifold $O$ as in Figure \ref{SXCTorOrb}. As we have previously seen, $O$ is a rigid orbifold. Therefore, $S'$ is totally geodesic and $L$ is a TGS link in $L(m,n)$.

Next, we would again like to extend $L$ to an infinite family of TGS links in $L(m,n)$. This can be achieved using the same three extension techniques as in Section \ref{SXC}. That is, we can add layers, add edges to the polygon in each layer, and add twists to each arm. However, these additions will have restrictions based on $L(m,n)$ and the values of $m$ and $n$. First, we will consider adding layers. Since we need to preserve the $\frac{1}{n}$ meridinal and $\frac{1}{m}$ longitudinal rotational symmetry $R$, we can add any multiple of $2m$ layers to $L$. Let $L'_\alpha$ be the link formed by adding $\alpha(2m)$ layers to $L$ and attaching the layers with arms as before. After this alteration, each layer of $L'_\alpha$ still has three arms. Therefore, the representativity of $L'_\alpha$ is at least 6 and $L'_\alpha$ is hyperbolic. Let $S'_\alpha$ be the surface comprised of the interior of each layer polygon and the interior of each arm. Then, $R$ is preserved and we can take the quotient of $S'_\alpha$ by it. This will yield $S'$ like in the original case of $L$. Therefore, since $S'$ is totally geodesic, it follows that $S'_\alpha$ is totally geodesic as well. 

Next, we can add edges to each layer polygon. Again, in order to preserve $R$, the number of edges added to each polygon must be a multiple of $2n$. Let $L''_\beta$ be the link obtained by adding $\beta(2n)$ edges to the layer polygon in each layer of $L$. Adding arms to each layer increases the lower bound on the representativity of $L''_\beta$. Therefore, $L''_\beta$ is hyperbolic. Let $S''_\beta$ be the surface formed by the interior of each layer polygon and each arm. Since $R$ is still a symmetry of $S''_\beta$, we can take the quotient by $R$ which will give us a pair of layer polygons connected by arms. This surface will also have rotational symmetries by lines which bisect edges of these polygons and symmetries which pierce the central crossings of arms. Hence, we can follow the same procedure as with $L$ in order to create a rigid orbifold. Hence, $S''_\beta$ is totally geodesic and $L''_\beta$ is a TGS link.

Third, we can add full-twists to each arm. Let $L'''_\gamma$ be the link we get by adding $\gamma$ full twists to each arm of $L$. In this case, the representativity of $L'''_\gamma$ is at least 6 since the quantity of arms has not changed. Thus, $L'''_\gamma$ is hyperbolic. Let $S'''_\gamma$ be the surface formed by the interior of each layer polygon and each arm. This alteration has no impact on $R$ or the symmetries about lines which bisect edges of the layer polygons. Also, by adding full-twists, we ensure that each arm will still have a crossing at its center point. Therefore, we retain all of the symmetries which pierce crossings of arms. Thus, we can use the same procedure as with $L$ to give us a rigid orbifold. Therefore, $S'''_\gamma$ is totally geodesic and $L'''_\gamma$ is a TGS link in $L(m,n)$.

Lastly, just like we did in Section \ref{SXC}, we can combine all three of these extention techniques and apply them simultaneously. Let $L_{\alpha \beta \gamma}$ be the link obtains by beginning with $L$ and then adding $\alpha(2m)$ layers, adding $\beta(2n)$ edges to the polygon in each layer, and adding $\gamma$ full-twists to each arm. By all of the arguments from above, $L_{\alpha \beta \gamma}$ is hyperbolic. Let $S_{\alpha \beta \gamma}$ be the surface made up of the interiors of each polygon and each arm. We saw above that each of these extensions preserves $R$ as well as the other symmetries used to quotient $S_{\alpha \beta \gamma}$ down to a rigid orbifold. Hence, every $S_{\alpha \beta \gamma}$ is totally geodesic and together all of the $L_{\alpha \beta \gamma}$ form an infinite family of TGS links in the lens space $L(m,n)$.

\subsection{Knotted Ribbon Links}

\begin{figure}[htbp]
\begin{center}
\includegraphics[width=5cm]{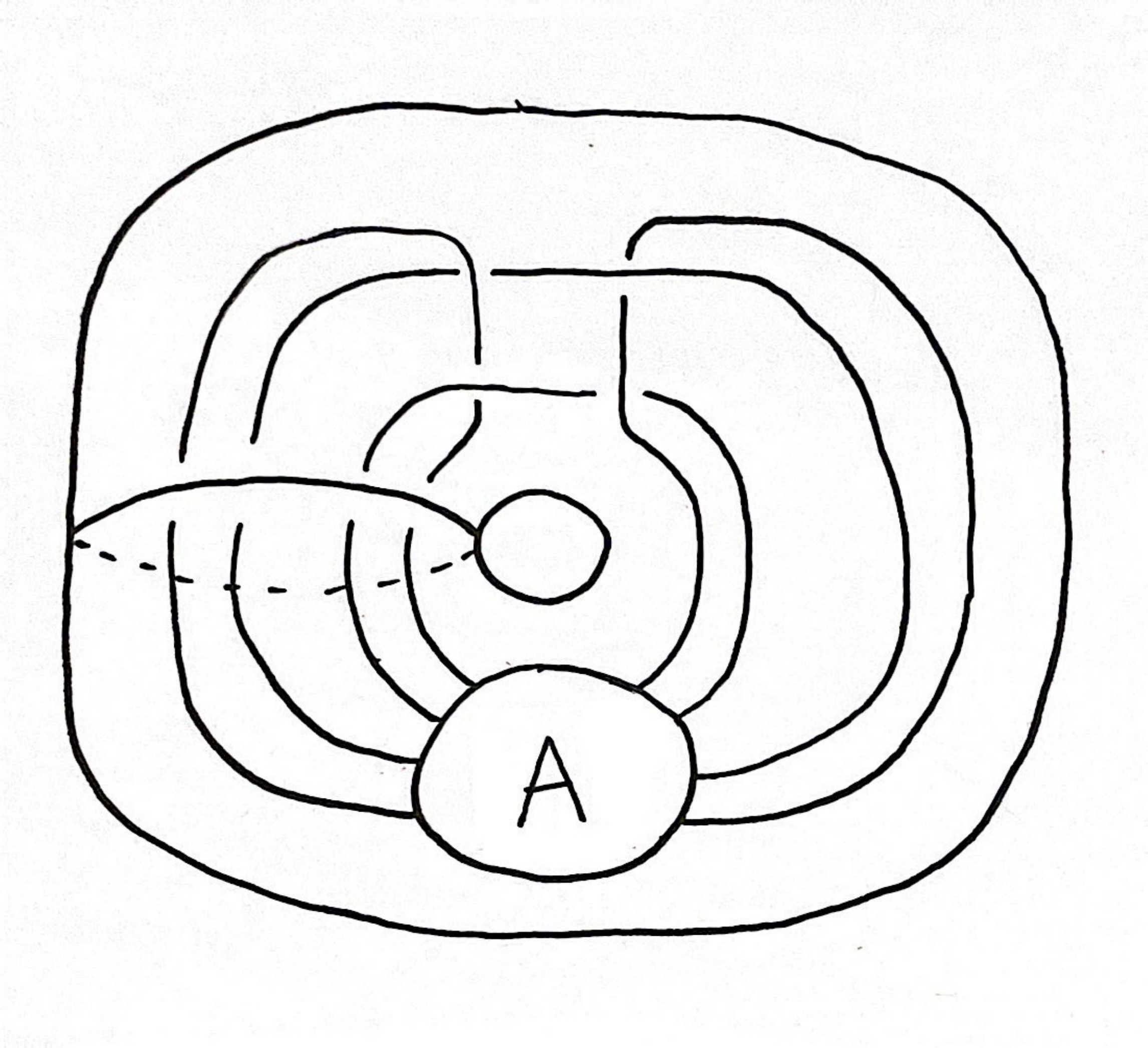}

\caption{}
\label{LensKR}
\end{center}
\end{figure}

We will now consider the second approach to finding totally geodesic spanning surfaces that does not utilize rigid orbifolds. As mentioned in \ref{intro}, Adams' work in \cite{thrice} tells us a surface embedded in a hyperbolic 3-manifold which is isotopic to a thrice-punctured sphere is totally geodesic. Therefore, in this subsection we will aim to construct hyperbolic links with spanning surfaces isotopic to a thrice-punctured sphere. We will construct a link $L$ inside a solid torus $\mathcal{T}$ that projects onto the boundary torus $T$. Then, an $(m,n)$ Dehn surgery on a meridian of $T$ yields the lens space $L(m,n)$ with $L$ embedded inside. Begin with a square at the top of $T$, and then pull the four pairs of strands at each crossing down to the bottom of $T$. We will call these strand pairs ribbons. Next, use the ribbons to form a tangle $A$, treating each ribbon as a link strand. Figure \ref{LensKR} shows the general structure of a link in this family.

At this time, we have not proven that these links are hyperbolic in $\mathcal{T}$. In fact, for many choices of tangle, the link will not be hyperbolic in $\mathcal{T}$. We could prove hyperbolicity by proving $\mathcal{T} \setminus L$ contains no essential disks, annuli, spheres, or tori. However, we are able to reasonably conjecture that a given link is hyperbolic through the use of Jeff Weeks' computer program SnapPy. After inserting a given link into SnapPy, the program is able to calculate the hyperbolic volume of the link. If this volume is a reasonable value, goes out to sufficiently many decimal places, and remains consistent upon multiple re-entries, then we are able to conjecture that the link is in fact hyperbolic. Table \ref{hypVols} shows a table with a number of these knotted ribbon links and the hyperbolic volume reported by SnapPy.

\begin{table}[htbp]
        \begin{tabular}{|c|c|c|c|}
            \hline
            Link & Hyperbolic Volume & Link & Hyperbolic Volume \\
            \hline
            \includegraphics[width=2.3cm]{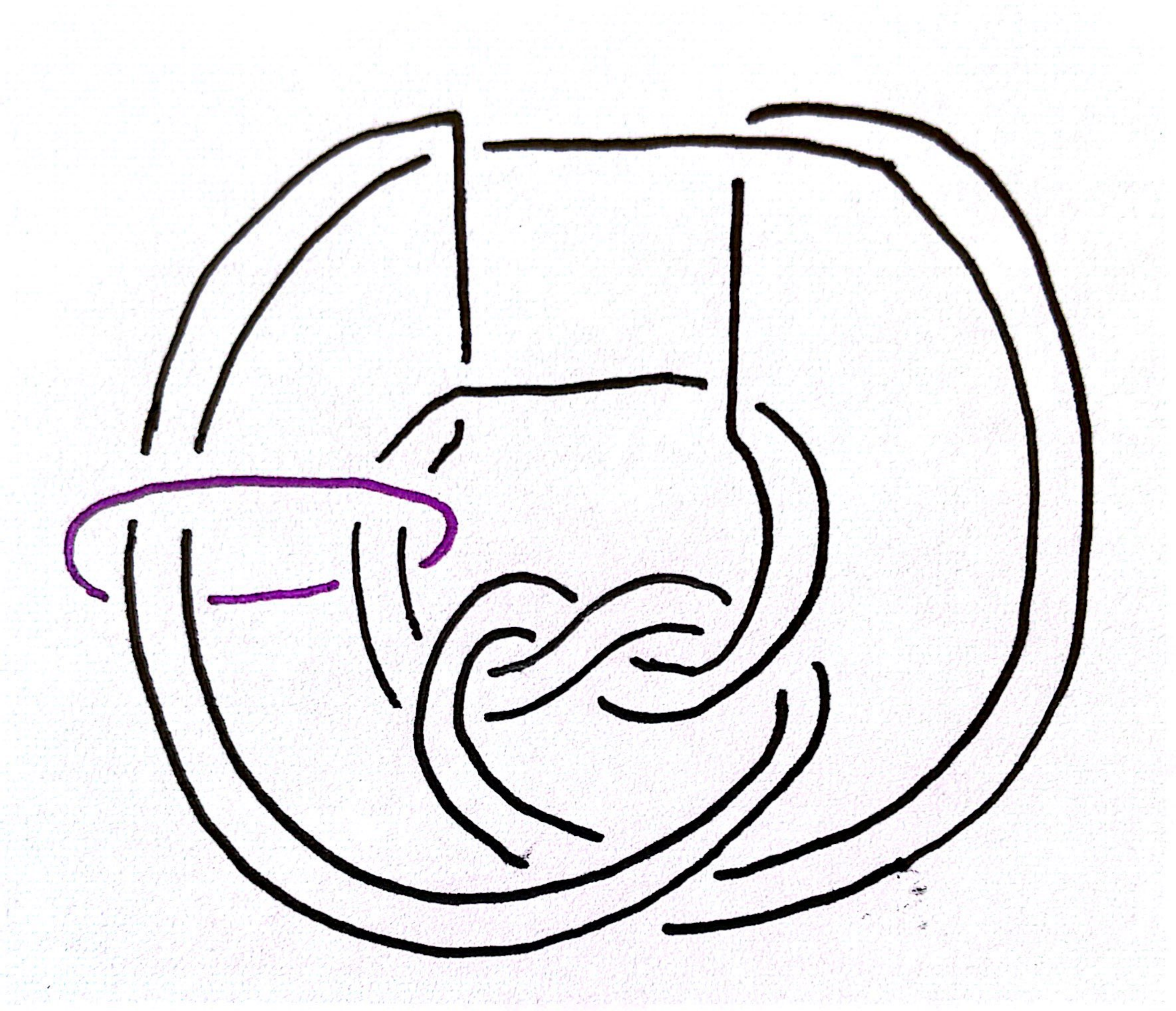} & 17.509452564 & \includegraphics[width=2.3cm]{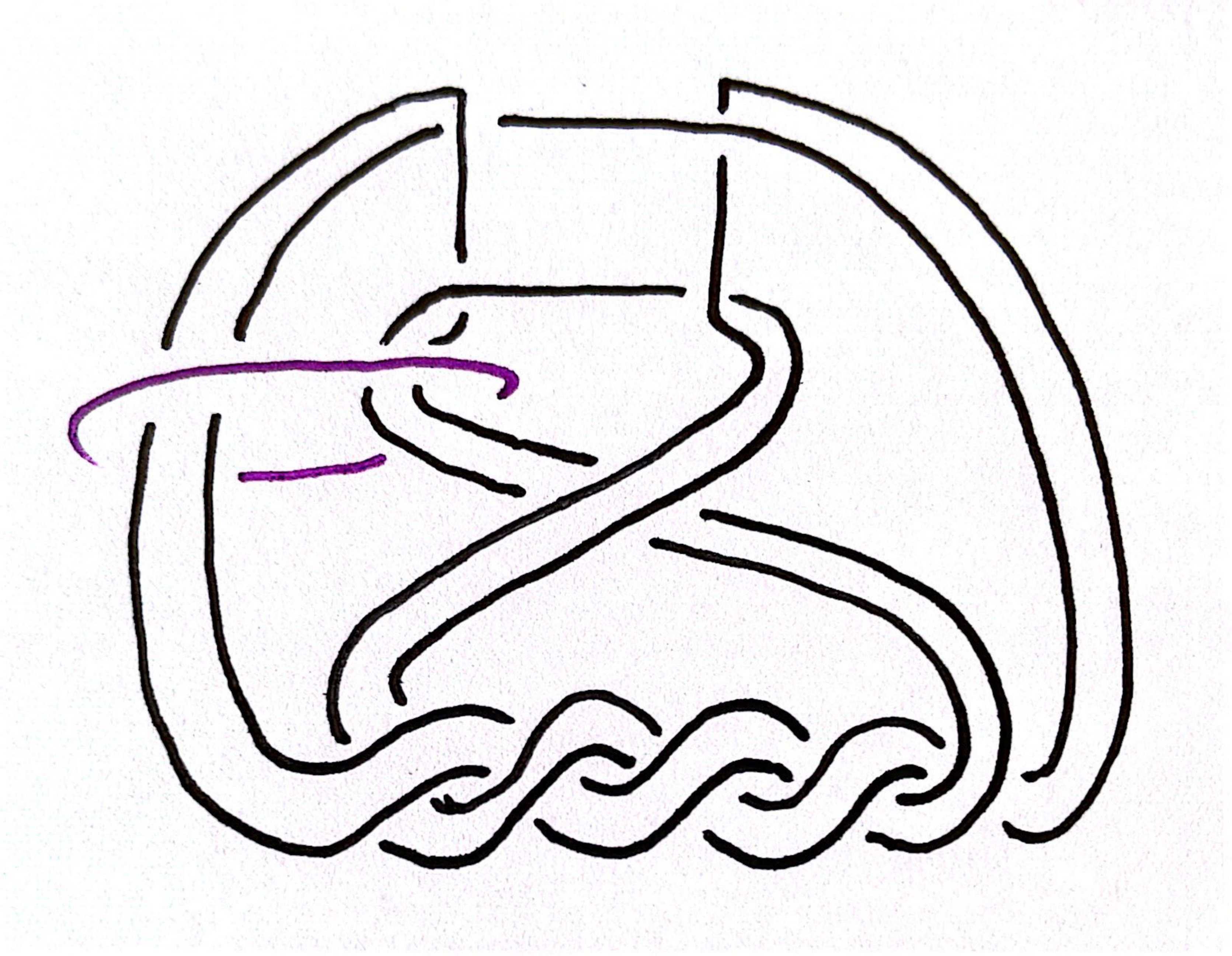} & 19.377635979 \\
            \hline
            \includegraphics[width=2.3cm]{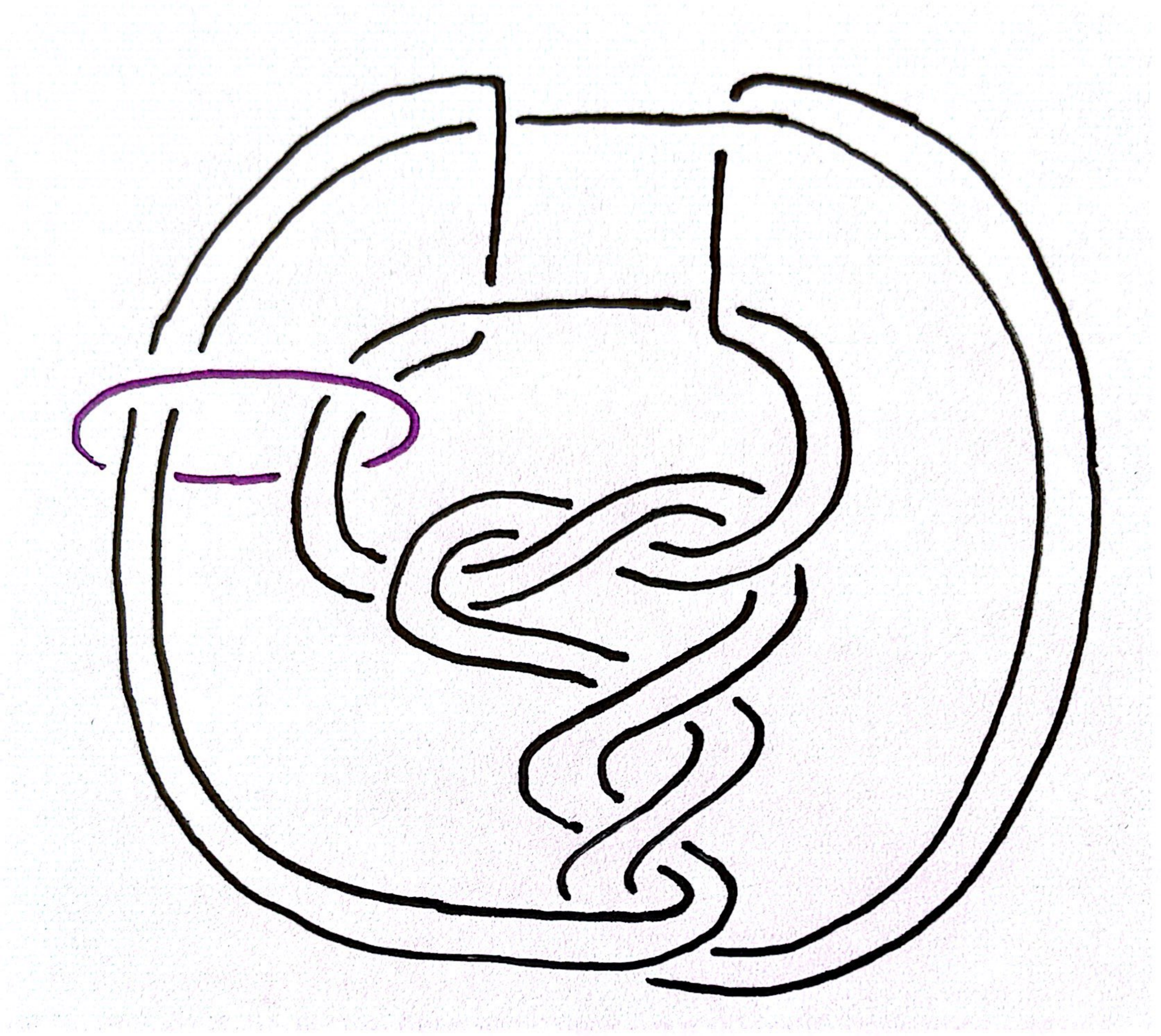} & 20.068375558 & \includegraphics[width=2.3cm]{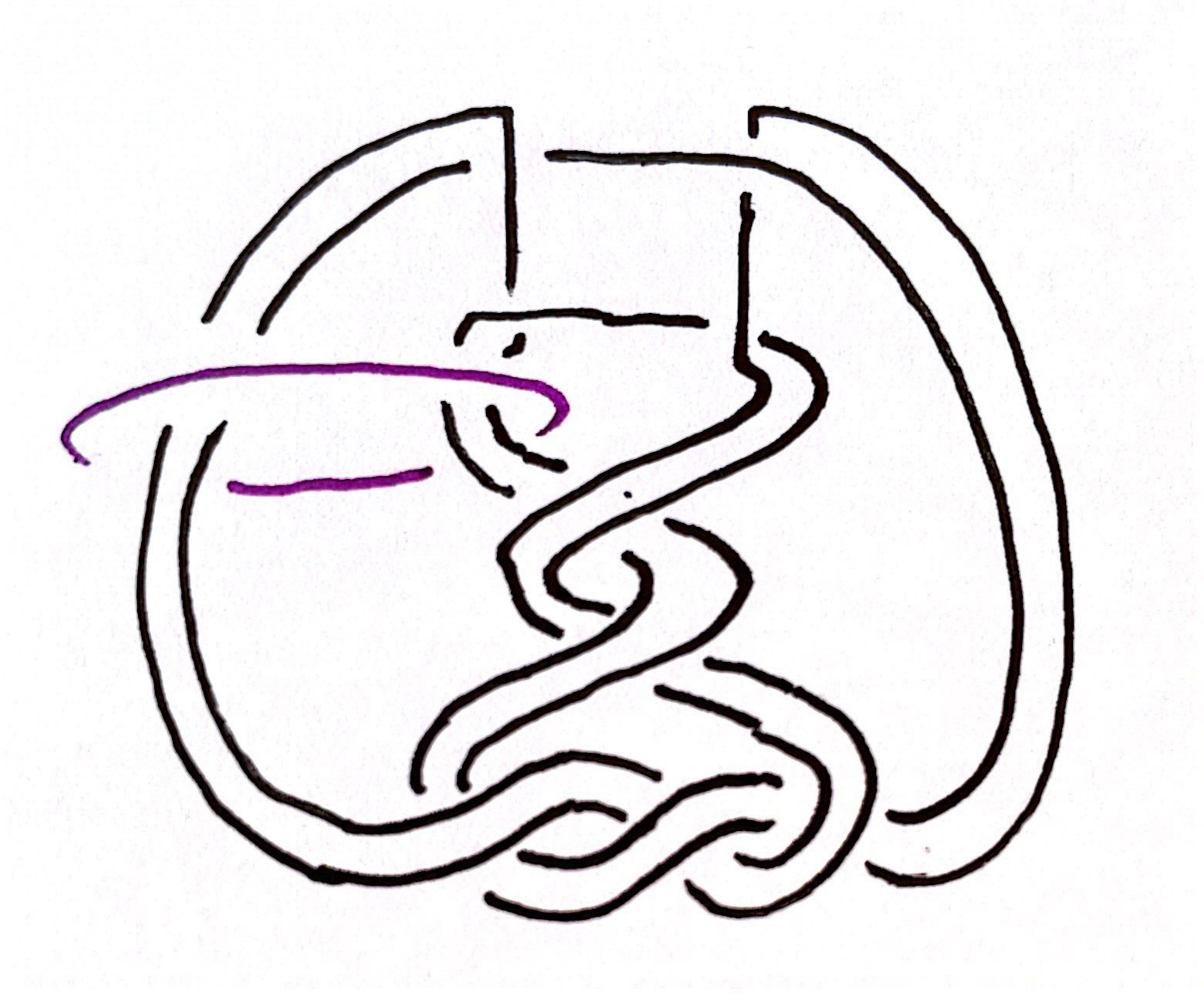} & 20.068375558 \\
            \hline
            \includegraphics[width=2.3cm]{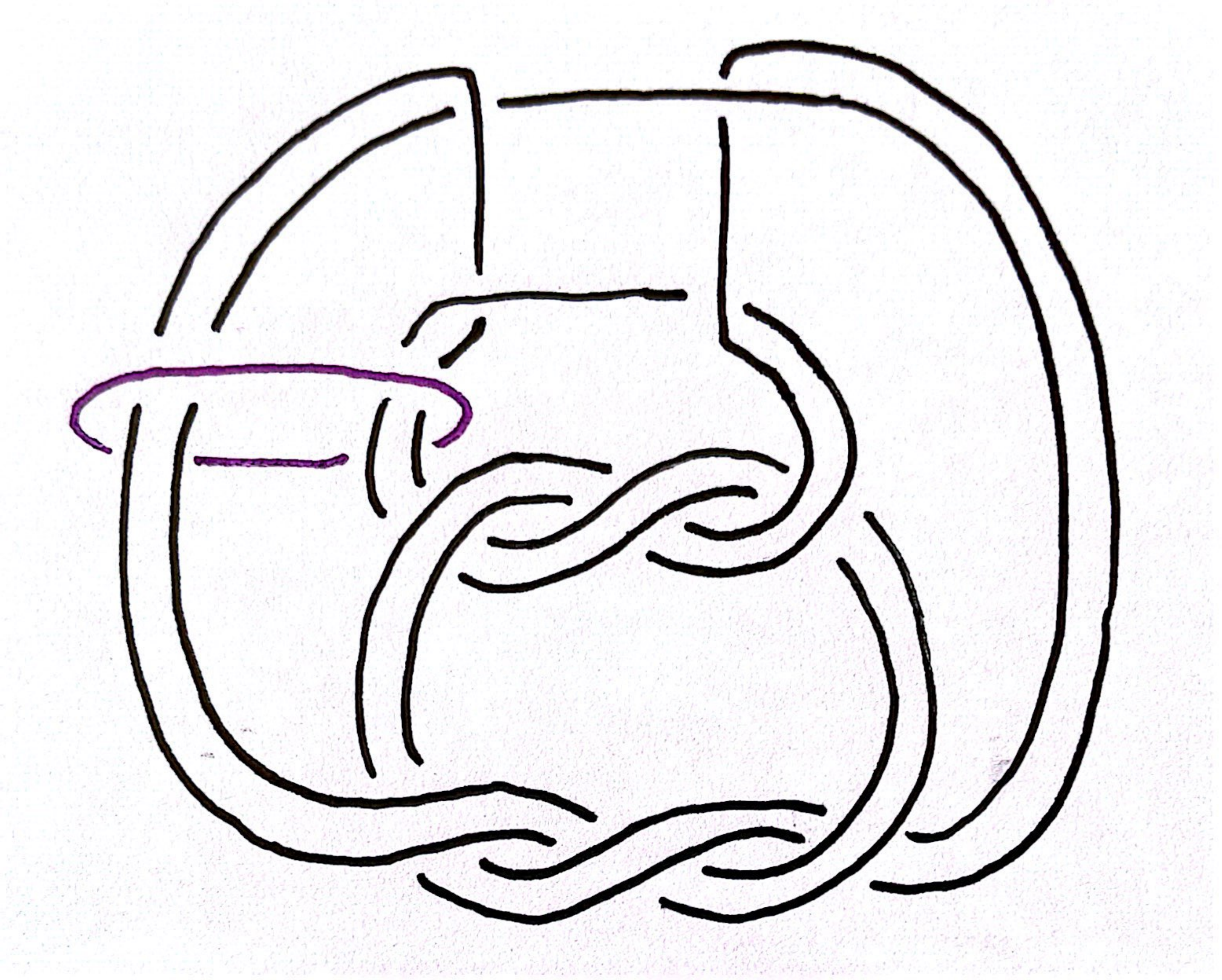} & 21.7380936554 & \includegraphics[width=2.3cm]{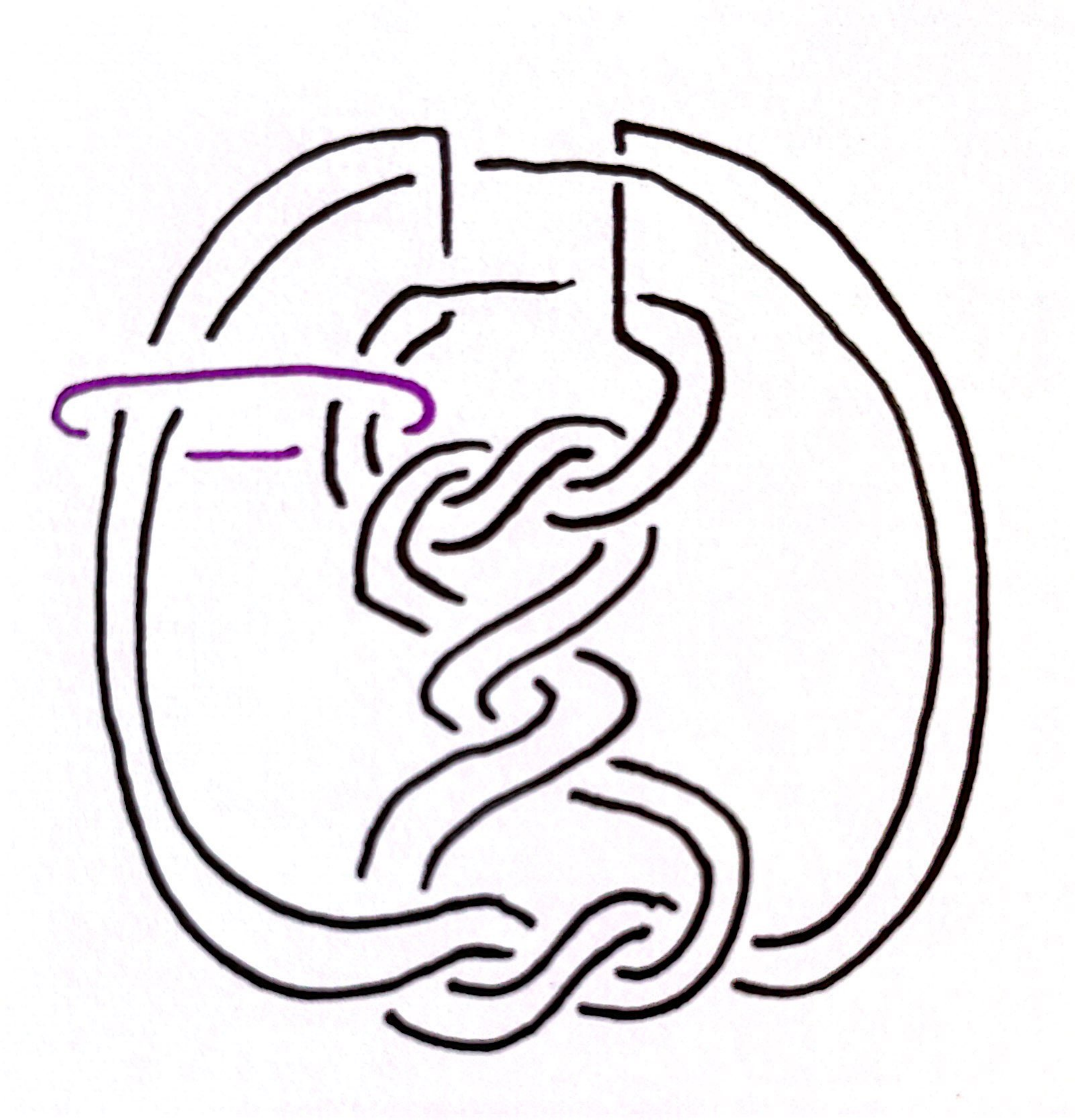} & 26.067864330 \\
            \hline
        \end{tabular}
        \caption{}
        \label{hypVols}
\end{table}

When depicted in $S^3$, all of these knotted ribbon links have four cusps. Three cusps form the link itself, and a fourth cusp which corresponds to a meridian of $T$ is added to project the link onto $T$. We will call this fourth cusp $C$, and each picture in Table \ref{hypVols} has $C$ drawn in purple. As we mentioned before, we can perform a $(m,n)$ Dehn surgery on $C$ in order to send $L$ to some link $L'$ in the lens space $L(m,n)$. Then, we can apply Theorem \ref{Thurston}. Since $L \cup C$ is hyperbolic in $S^3$, the result of performing a Dehn surgery on $C$ will yield a hyperbolic link in all but finitely many cases. Hence, $L'$ will be hyperbolic in all but finitely many lens spaces. 

Next, given a knotted ribbon link $L$, let $S$ be the spanning surface of $L$ formed by the interior of the top square and the interiors of the two ribbons. Recall, we designed $L$ so that $S$ would topologically be a thrice-punctured sphere, and this is the case. Therefore, $S$ is totally geodesic in $\mathcal{T}$ by \cite{thrice} Then, for all of the Dehn surgeries which yield a hyperbolic link $L'$ in the associated lens space, $S$ will lift to a totally geodesic surface $S'$. Therefore, $L'$ is a TGS link in its corresponding lens space. Hence, we get infinitely many knotted ribbon links which yield TGS links in infinitely many lens spaces.

We can again use SnapPy to see their totally geodesic structure explicitly by examining the horoball diagram associated with these links $L'$ which have undergone Dehn surgery. In a horoball diagram, totally geodesic surfaces correspond to straight lines of full-size horoballs. We will look at the first link in Table \ref{hypVols}, call it $X$, and its associated link in a number of different lens spaces. 

\begin{figure}
    \centering
    \begin{tabular}{cc}
       \includegraphics[width=5.5cm]{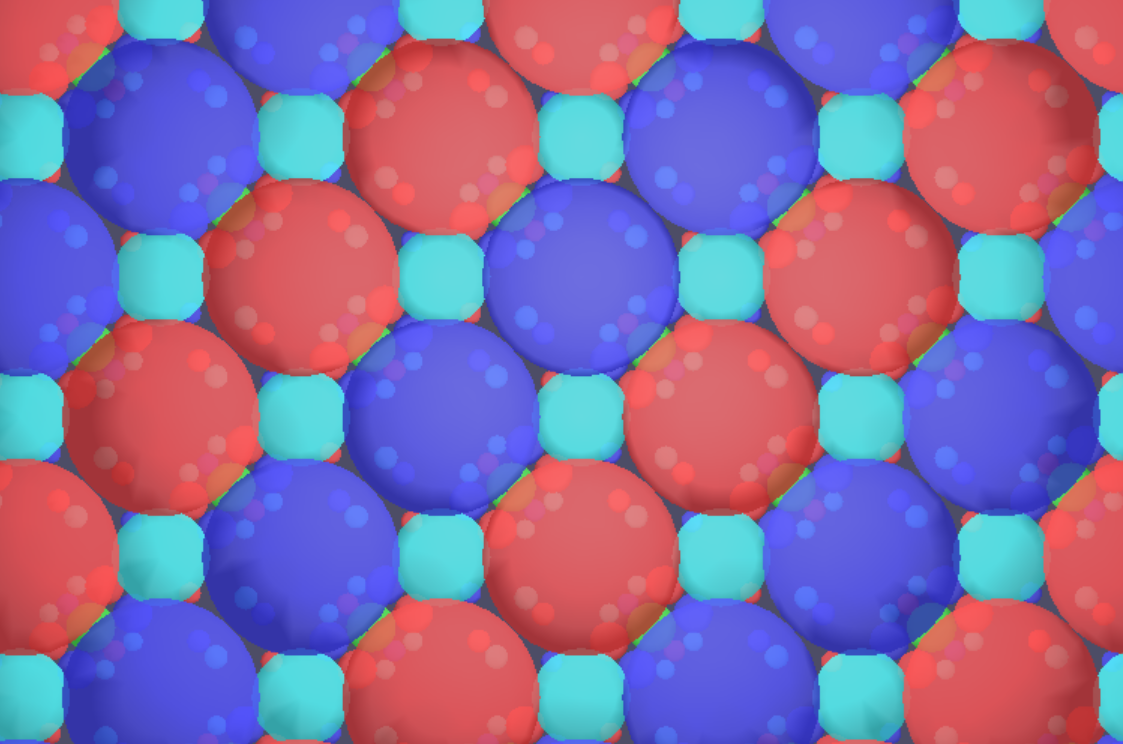}  & \includegraphics[width=5.5cm]{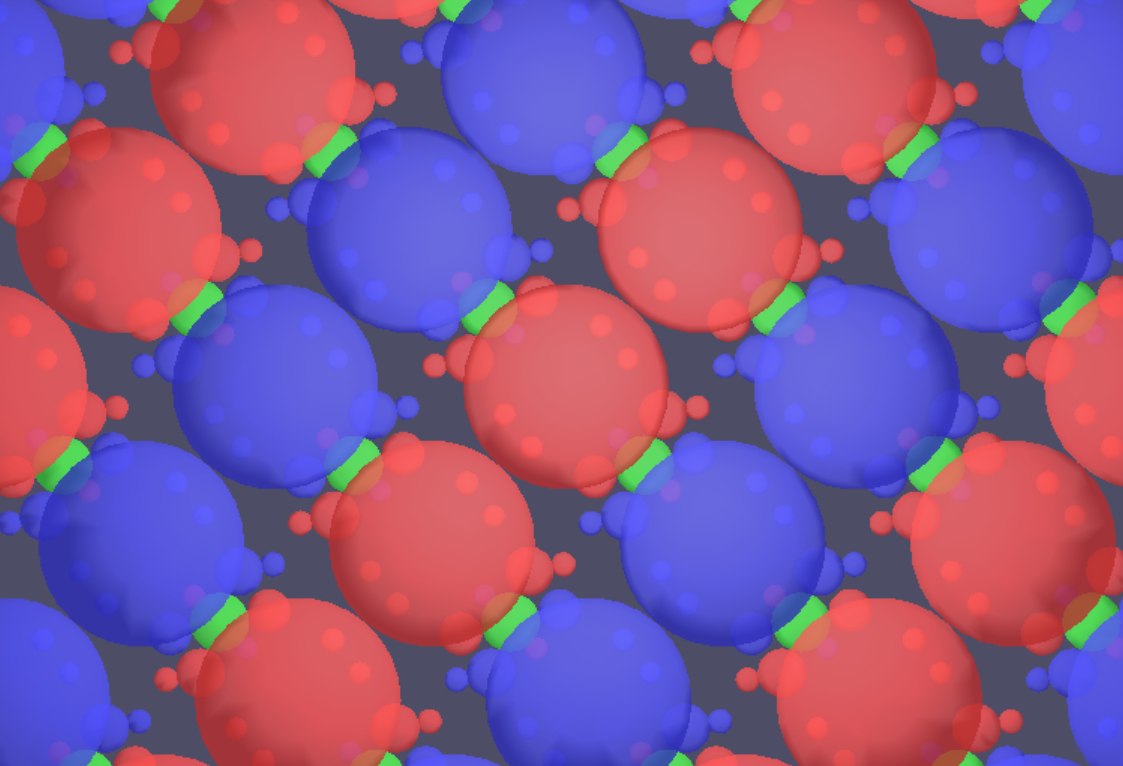} \\
       \includegraphics[width=5.5cm]{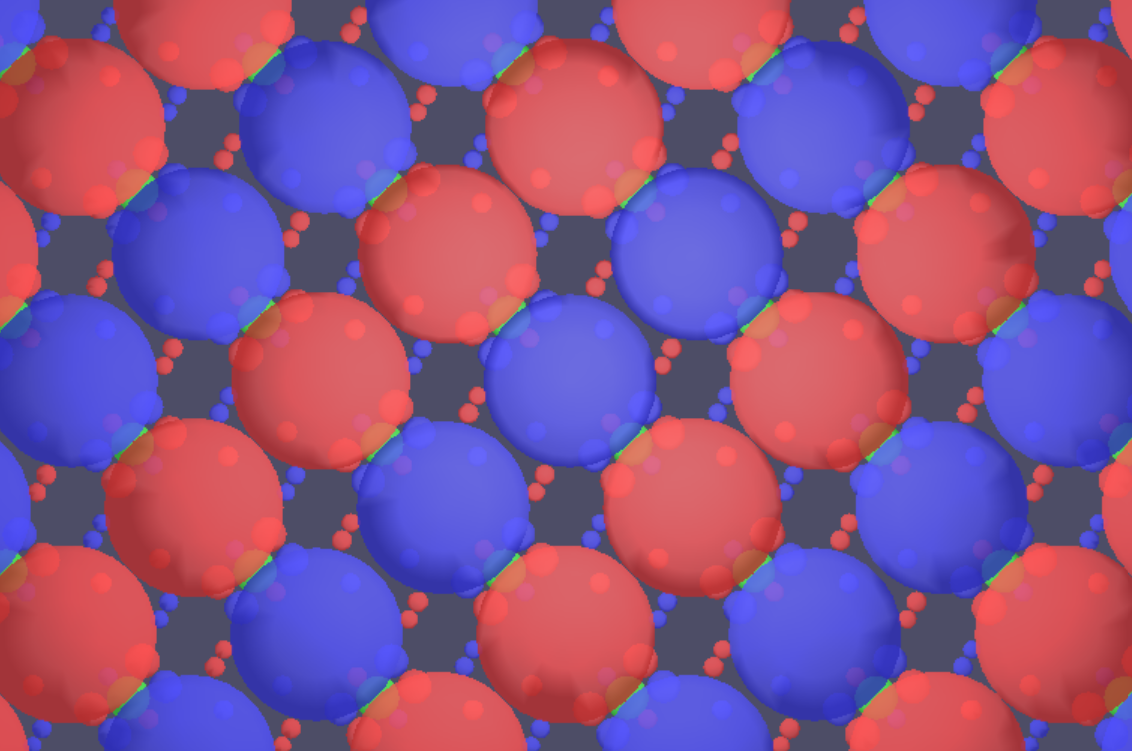}  & \includegraphics[width=5.5cm]{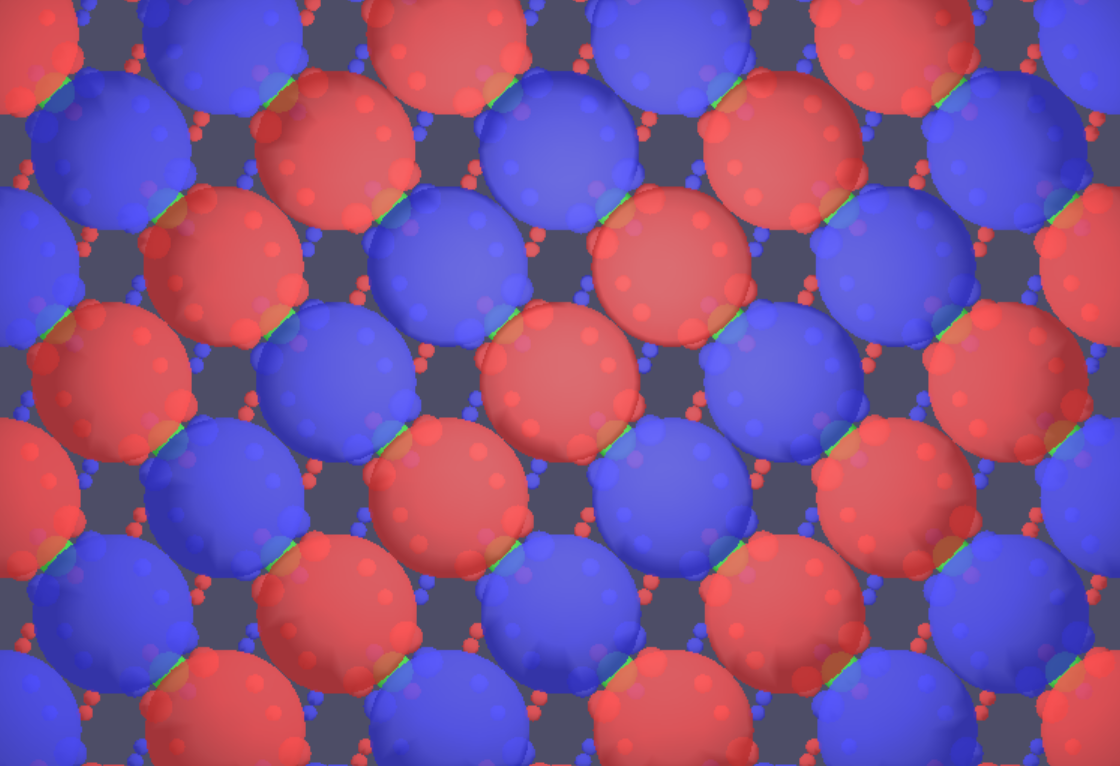}
    \end{tabular}
    
    \caption{}
    \label{LensHB}
\end{figure}

Figure \ref{LensHB} shows the horoball diagrams of $X$ before any Dehn surgery and $X'$ after $(1,3), (5,2)$, and $(3,4)$ Dehn surgeries. The aqua horoballs correspond to the cusp labeled $C$ before, on which we performed the surgeries. These diagrams all include multiple straight lines of horoballs, but it is most likely that the diagonal lines running from the top left to the bottom right that include blue, red, and green horoball corresponds to our surface $S$ and the associated $S'$ in each lens space.

\section{TGS Links with Non-Orientable Spanning Surfaces}

\begin{figure}[htbp]
\begin{center}
\includegraphics[width=7cm]{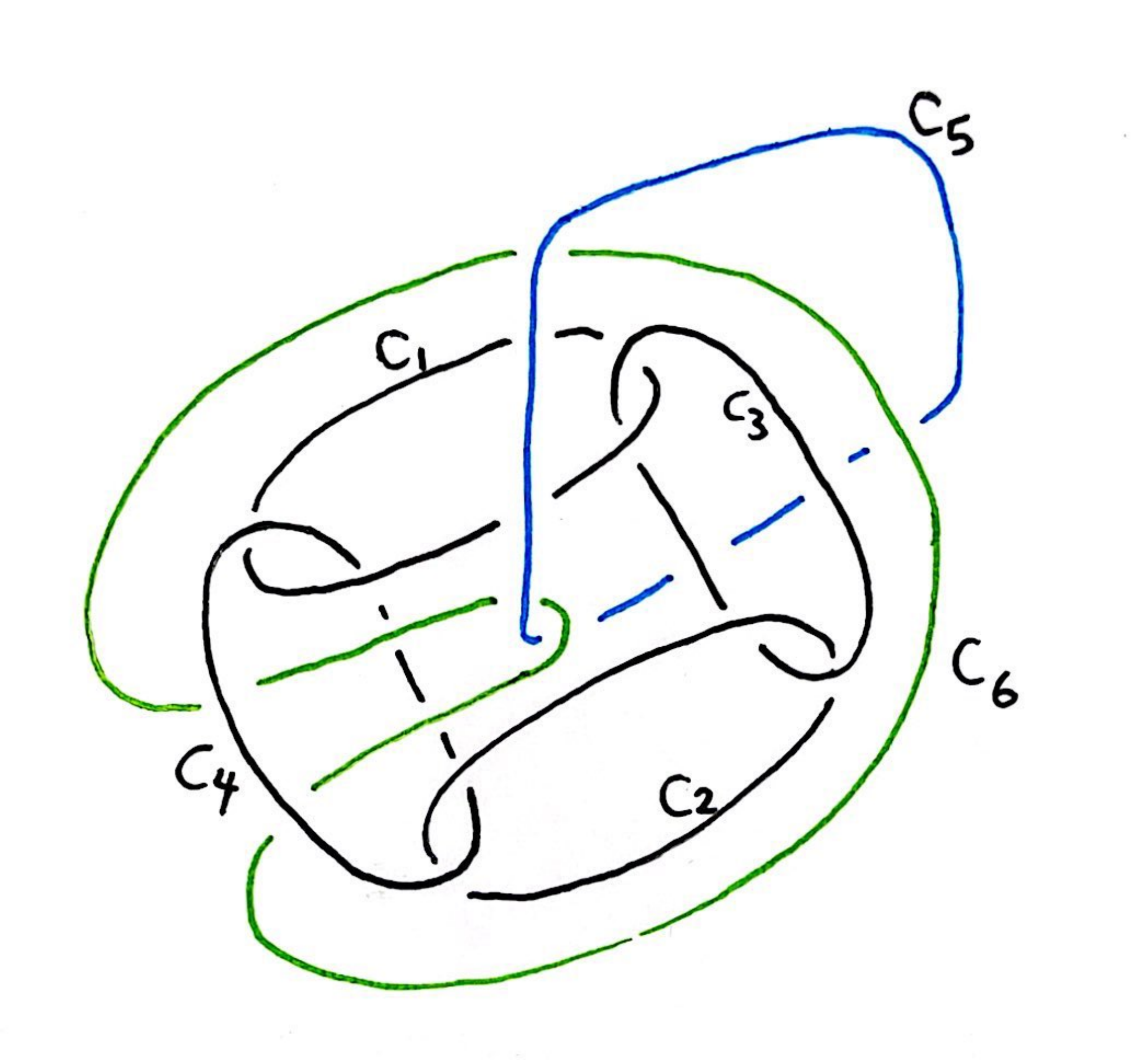}

\caption{}
\label{NonOr}
\end{center}
\end{figure}

Next, we will look at links in the 3-sphere $S^3$ and the solid torus $\mathcal{T}$ in order to find non-orientable surfaces which are totally geodesic. As we have in previous sections, we will build a link in $\mathcal{T}$ by beginning with a link in $S^3$ and then adding a trivial link component whose complement is a torus. A depiction of this link can be seen in Figure \ref{NonOr}.

To create the link $L$, we will start with 2 trivial links components $C_1, C_2$ that lie in some plane $P$. Then, we will add 2 more trivial link components $C_3, C_4$ each of which are perpendicular to $P$ and run through both $C_1$ and $C_2$. Note, this is the minimally twisted chain of 4 components in $S^3$. Next, we will add a trivial link component $C_5$, drawn in blue, which will be perpendicular to $P$ and run through the interior and exterior of the chain. Lastly, we will add a final trivial link component $C_6$, drawn in green, that lies in $P$. We will let $C_6$ surround $C_1, C_2,$ and $C_3$, and we will pass through the middle of $C_4$ twice. Note, $C_5$ and $C_6$ are not linked with one another. The link components $C_3, C_4,$ and $C_5$ are all in fact perpendicular to $P$ since reflecting $L$ through $P$ sends each of these components back to themselves.

\begin{figure}[htbp]
\begin{center}
\includegraphics[width=10cm]{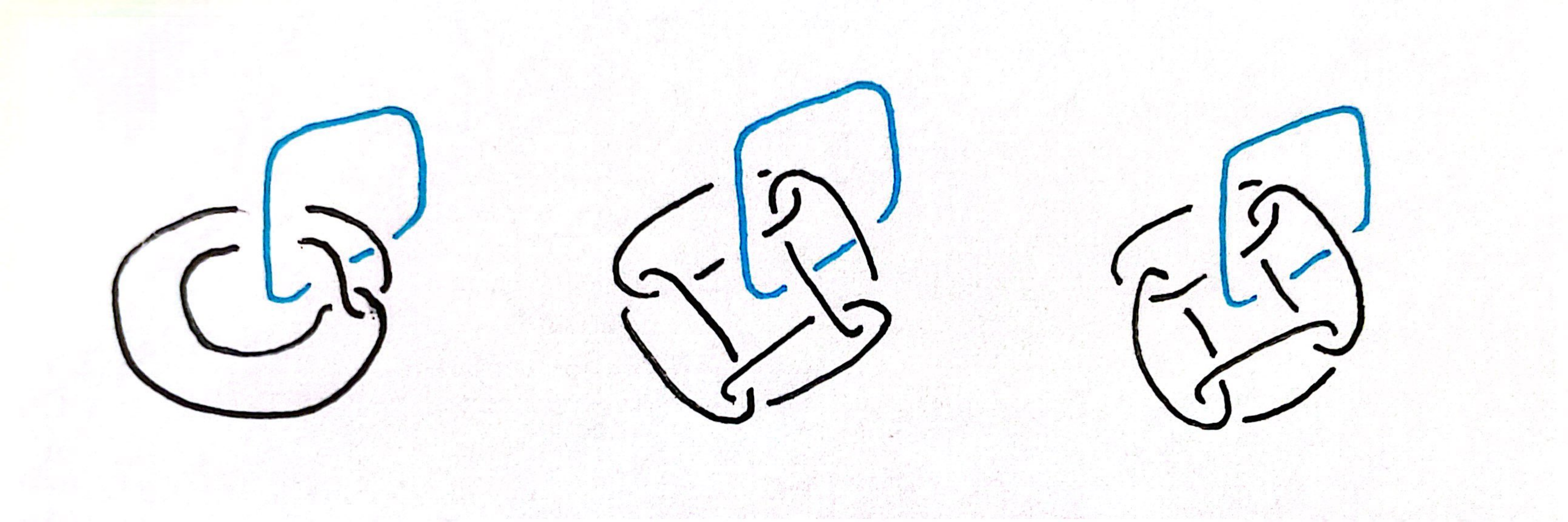}

\caption{}
\label{NonOrChains}
\end{center}
\end{figure}

Now, we will show $L$ is hyperbolic. Begin with the Whitehead link $W$, shown in the first diagram of Figure \ref{NonOrChains} This link is known to be hyperbolic. Then, take a 4-fold cover of $W$ about the blue component $B$. This new link can be seen in the second diagram of Figure \ref{NonOrChains}, and call it $W'$. Taking a finite cover of a hyperbolic link yields a hyperbolic link, so $W'$ is hyperbolic. Next, consider the surface $S$ bounded by $B$. This surface $S$ is a twice-punctured disk, so it is totally geodesic by Theorem \ref{AdamsThrice}. Therefore, we can cut along $S$, perform a half twist, and glue $S$ back to itself two times to change $W'$ into the minimally twisted chain of four components, seen in the third diagram of Figure \ref{NonOrChains}. Call this link $W''$. Note, $W''$ is also hyperbolic since taking a hyperbolic knot and cutting, twisting, and gluing along a thrice-punctured sphere yields a hyperbolic link. Lastly, we will add one more component to $W''$, $C_6$ as defined above. This gives us our link $L$. Consider $M = S^3 \setminus W''$, which is an orientable, hyperbolic 3-manifold with finite volume. Consider the surfaces in $P$ to the inside of $C_6$ and to the outside of it. The inside surface is a twice-punctured disk, with $C_1$ and $C_3$ forming those punctures. The outside surface is also a twice-punctured disk, with $C_5$ forming both punctures. Therefore, $C_6$ is a curve which separates the punctures of a four-punctured sphere into two pairs. Hence, $C_6$ is an embedded geodesic in $M$. Therefore, $M \setminus C_6$ admits a hyperbolic metric by Theorem \ref{KojimaGeodisic}. Therefore, $L$ is hyperbolic.

\begin{figure}[htbp]
\begin{center}
\includegraphics[width=5cm]{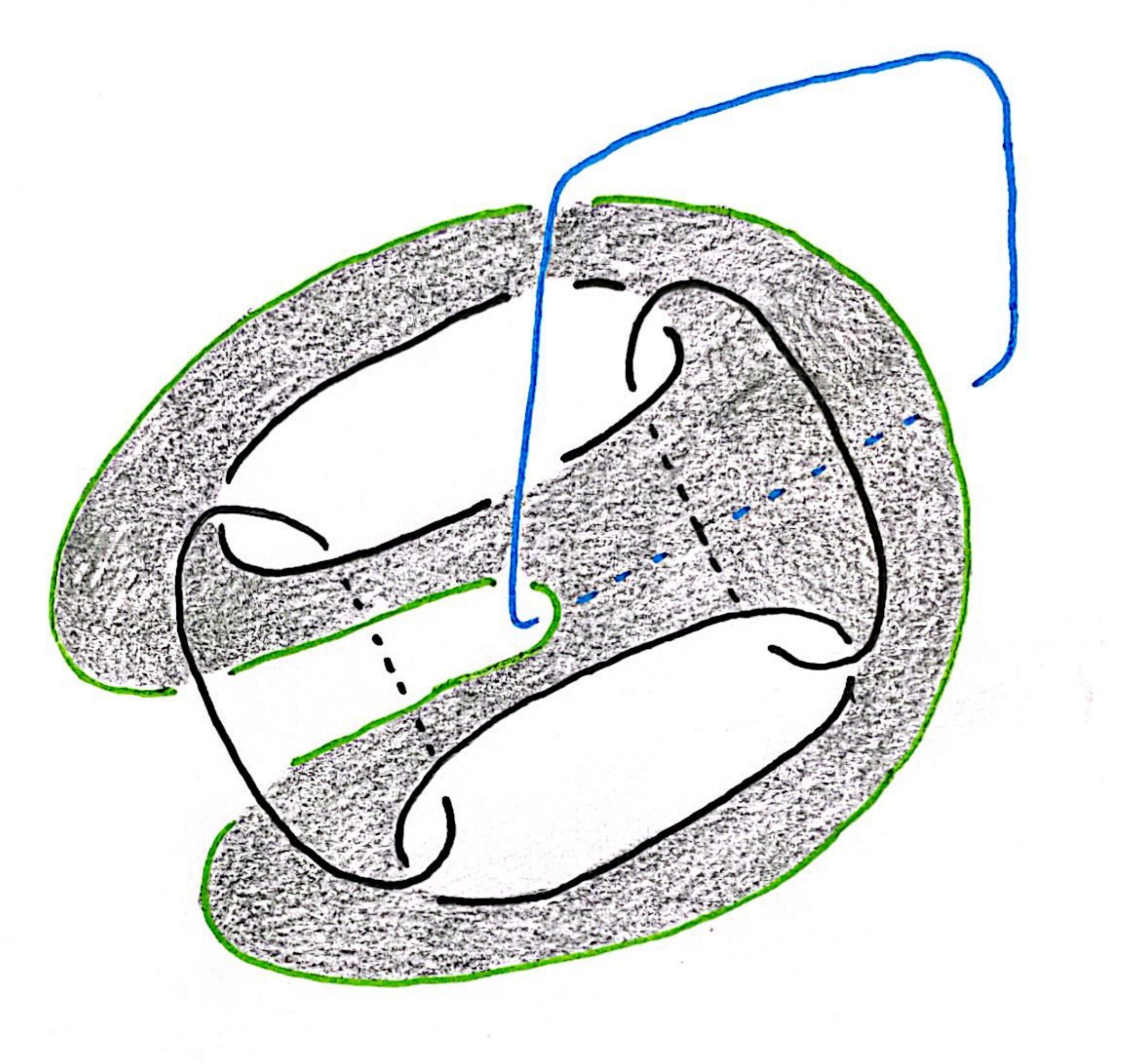}

\caption{}
\label{NonOrSurf}
\end{center}
\end{figure}

Let $S$ be the surface which lies in $P$ and is the interior of $C_6$ and punctured by $C_1$ and $C_3$, as in Figure \ref{NonOrSurf}. Since $S$ has boundary along 3 cusps, $S$ is topologically a twice-punctured disk or equivalently a thrice-punctured sphere. Therefore, $S$ is totally geodesic by Theorem 1.5. Thus, $L$ is a TGS link in $S^3$ and the link formed by $C_1, C_2, C_3, C_4, C_6$ is a TGS link in $\mathcal{T}$. 

\begin{figure}[htbp]
\begin{center}
\includegraphics[width=5cm]{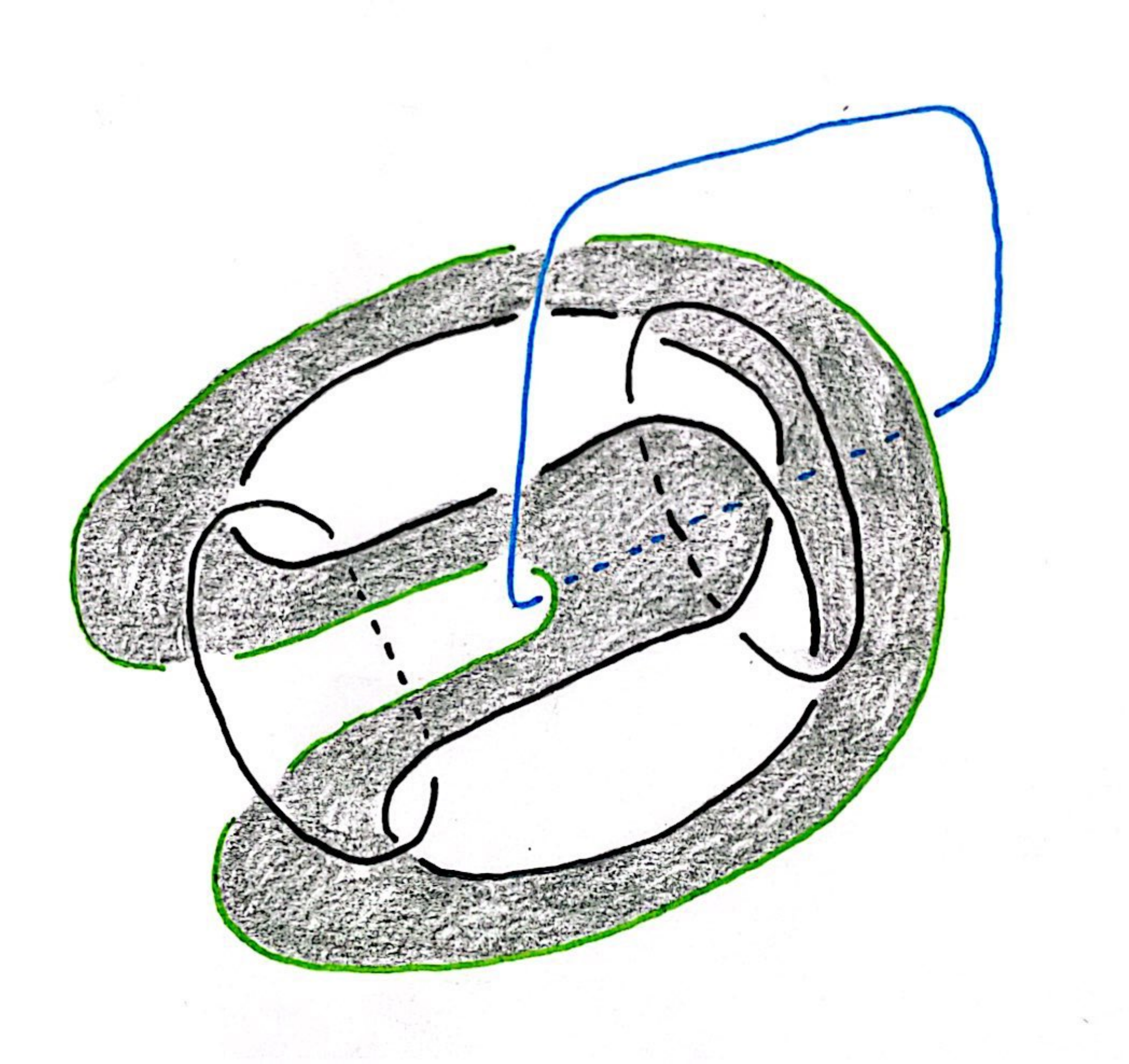}

\caption{}
\label{NonOrSurfTw}
\end{center}
\end{figure}

Next, we will transform $L$ so that $S$ becomes a non-orientable surface while remaining totally geodesic. Note, the surface $F$ formed by the interior of $C_2$ is disk which is punctured a total of two times by $C_1$ and $C_3$. Hence, $F$ is topologically a thrice-punctured sphere and thus is totally geodesic by Theorem \ref{AdamsThrice}. Additionally, since $C_2$ is perpendicular to the projection plane $P$, we are able to cut along $F$ and glue $F$ back to itself while preserving the geometry of the rest of the manifold. Hence, we will cut along $F$, rotate $F$ $180^\circ$, and glue it back to itself to create a new link $L'$ where $C_1$ and $C_3$ have merged into a single component with a half-twist. Let $S'$ be the transformed $S$. Both $L'$ and $S'$ can be seen in Figure \ref{NonOrSurfTw}. Now, $S'$ is non-orientable and is also totally geodesic since the transformation involved cutting along a totally geodesic surface perpendicular to $S'$. Hence, $L'$ is a TGS link in $S^3$ torus with a non-orientable spanning surface, and $L' \setminus C_5$ is a TGS link in $\mathcal{T}$ with a non-orientable spanning surface.

\begin{figure}[htbp]
\begin{center}
\includegraphics[width=5.5cm]{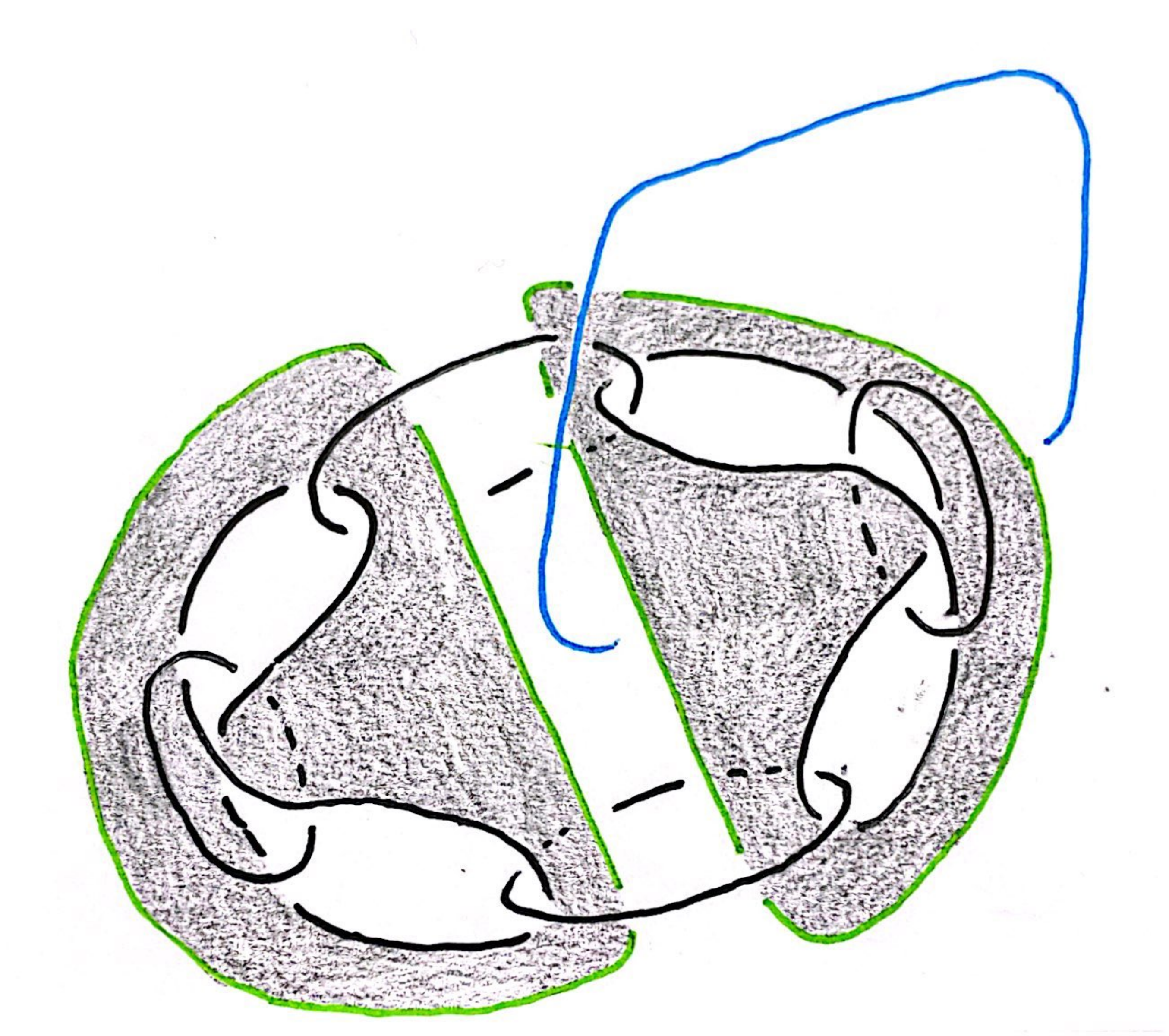}

\caption{}
\label{NonOrSurf2}
\end{center}
\end{figure}

We can extend $L$ to a family of links in $\mathcal{T}$ which can be transformed into TGS links by taking a double cover over $C_5$. Let $L'_2$ be the link obtained by taking a double cover of $L'$ over $C_5$. Note, $L'_2$ is hyperbolic since the double cover of a hyperbolic link is hyperbolic. This yields two non-orientable totally geodesic spanning surfaces, which can be seen in the Figure \ref{NonOrSurf2}. Furthermore, we can take an $n$-fold cover of $L'$ a link $L'_n$ which will have $n+1$ non-orientable totally geodesic spanning surfaces.

\bibliographystyle{plain}
\bibliography{references}

\end{document}